\input{aipcheck}

\documentclass[
  ,final            
  ]
  {aipproc}

\layoutstyle{6x9}

\usepackage{longtable}
\usepackage{amssymb}
\usepackage{latexsym}
\usepackage{graphicx}
\usepackage{multicol}
\usepackage{multirow}
\usepackage{color}
\usepackage{pstricks,pst-plot,pst-text,pst-tree,pst-eps,pst-fill,pst-node,pst-math}

\newcommand{\dbl}{[\![}
\newcommand{\dbr}{]\!]}
\newcounter{numerothm}
\newcounter{numerodef}
\newcounter{numerostat}
\newcounter{numeroprop}
\newcounter{numerocor}

\newcommand{\proof}[1]{\textbf{\emph{Proof}} ~~ #1\\}
\newcommand{\thm}[1]{\stepcounter{numerothm} \textbf{{\emph{Theorem}} \thenumerothm} ~~#1\\}
\newcommand{\defn}[1]{\stepcounter{numerodef} \textbf{{\emph{Definition}} \thenumerodef} ~~#1\\}
\newcommand{\stat}[1]{\stepcounter{numerostat} \textbf{{\emph{Statement}} \thenumerostat} ~~#1\\}
\newcommand{\prop}[1]{\stepcounter{numeroprop} \textbf{{\emph{Proposition}} \thenumeroprop} ~~#1\\}
\newcommand{\cor}[1]{\stepcounter{numerocor} \textbf{{\emph{Corollary}} \thenumerocor} ~~#1\\}
\newcommand{\rem}[1]{\textbf{\emph{Remark}} ~~#1\\}
\newcommand{\nota}[1]{\textbf{\emph{Notation}} ~~#1\\}

\begin{document}

\title{Global Orbit Patterns for Dynamical Systems On Finite Sets }



\author{Ren\'e LOZI $^{\dag}$}{
  address={Laboratoire J.A. Dieudonn\'e - UMR du CNRS n$^{o}$6621\\
Universit\'e de Nice Sophia-Antipolis \\ Parc Valrose \\  06108 NICE
CEDEX 2 FRANCE\\},email=rlozi@unice.fr }

\author{Clarisse FIOL}{
  address={IUFM C\'elestin Freinet - Universit\'e de Nice Sophia-Antipolis\\ 89, av. George V\\
06046 NICE CEDEX 1 FRANCE \\},email=fiol@unice.fr  }

\begin{abstract}
In this paper, the study of the global orbit pattern (gop) formed by
all the periodic orbits of discrete dynamical systems on a finite
set $X$ allows us to describe precisely the behaviour of such
systems. We can predict by means of closed formulas, the number of
gop of the set of all the function from $X$ to itself. We also
explore, using the brute force of computers, some subsets of locally
rigid functions on $X$, for which interesting patterns of periodic
orbits are found.
\end{abstract}

\maketitle


\section{Introduction}\label{section1}

In some engineering applications such as chaotic encryption, chaotic
maps have to exhibit required statistical and spectral properties
close to those of random signals. There is a growing industrial
interest to consider and study thoroughly the property of such map
\cite{Patent2007, Patent2006, Patent2001}.

Very often, dynamical systems in several dimensions are obtained
coupling 1-dimension ones and their properties are strongly linked
\cite{Lozi1}.

Quasi-periodic or chaotic motion is frequently present in
complicated dynamical systems whereas simple dynamical systems often
involve only periodic motion. The most famous theorem in this field
of research is the Sharkovskii's theorem, which addresses the
existence of periodic orbits of continuous maps of the real line
into itself. This theorem was once proved toward the year 1962 and
published only two years after \cite{Shark1}.

Mathematical results concerning periodic orbits are often obtained
for functions on  real intervals. However, most of the time, as the
complex behaviour of chaotic dynamical systems is not explicitly
tractable, mathematicians have recourse to computer simulations. The
main question which arises then is: does these numerical
computations are reliable ?

As an example we report the results of some computer experiments on
the orbit structure of the discrete maps on a finite set which arise
when the logistic map is iterated "naively" on the computer. \\

Due to the discrete nature of floating points used by computers,
there is a huge gap between these results and the theoretical
results obtained when this map is considered on a real interval.
This gap can be narrowed in some sense (i.e. avoiding the collapse
of periodic orbits) in higher dimensions when ultra weak coupling is
used \cite{Lozi2, Lozi3}.

Nowadays the claim is to understand precisely which periodic orbit
can be observed numerically in such systems. In a first attempt we
study in this paper the orbits generated by the iterations of a
one-dimensional system on a finite set $X_N$ with a cardinal $N$.
The final goal of a good understanding of the actual behaviour of
dynamical systems acting on floating numbers (i.e. the numbers used
by computers) will be only reached after this first step will be
achieved.

On finite set, only periodic orbits can exist. For a given function
we can compute all the orbits, all together they form a global orbit
pattern. We formalise such a gop as the ordered set of periods when
the initial value thumbs the finite set in the increasing order. We
are able to predict, using closed formulas, the number of gop for
the set $\mathcal{F}_N$ of all the functions on $X$. We also explore
by computer experiments special subsets of $\mathcal{F}_N$, such as
sets of locally "rigid" functions which presents interesting
patterns of gop.

This article is organized as follows : in the section "Computational
divergences" we display some examples of such computational
divergences for the logistic map in various ways of discretization.
In the section "Pattern defined by all the orbits of a dynamical
system" we introduce a new mathematical tool: the global orbit
pattern, in order to describe more precisely the behaviour of
dynamical systems on finite sets. In the section "Cardinal of
classes" we give some closed formulas related to the cardinal of
classes of gop of $\mathcal{F}_N$. In the section "Functions with
local properties" we study the case of sets of functions with a kind
of local "rigidity" versus their gop, in order to show the
usefulness of these new tools.


\section{Computational divergences}\label{section2}

\subsection{Discretized logistic map}

\begin{figure}
\includegraphics*[width=6cm]{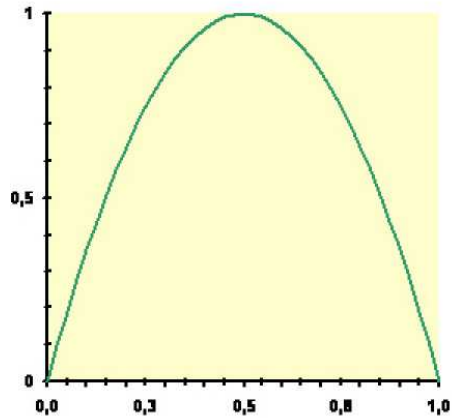}
\caption{Graph of the map $f(x)=4x(1-x)$ on $[0,1]$} \label{fig1}
\end{figure}

As an example of collapsing effects which happen when using
computers in numerical experiments, we presents the results of a
sampling study in double precision of a discretization of the
logistic map $f_4 : [0,1] \rightarrow [0,1]$ (see Fig. \ref{fig1})

\begin{equation}\label{eq1}
   f_4(x)=4x(1-x)
\end{equation}

and its associated dynamical system

\begin{equation}\label{eq2}
    x_{n+1}=4 x_n(1-x_n)
\end{equation}

which has excellent ergodic properties on the real interval.

\begin{figure}
\includegraphics*[width=6cm]{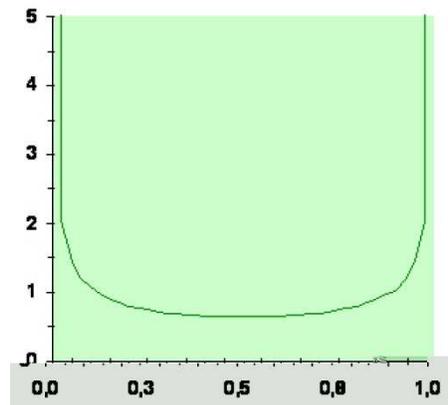}
\caption{Invariant measure of the logistic map} \label{fig2}
\end{figure}

There exists an unstable fixed point 0. \\
The set
$\left\{\frac{5-\sqrt{5}}{8},\frac{5+\sqrt{5}}{8}\right\}=\{0.3454915,0.9045084\}$
 is the period-2 orbit. \\
 In fact there exist an infinity of periodic
orbits and an infinity of periods. This dynamical system possesses
an invariant measure (see Fig. \ref{fig2}):

\begin{equation}\label{eq3}
   P(x)=\frac{1}{\pi \sqrt{x(1-x)}}
\end{equation}

However, in numerical computations using ordinary (IEEE-754) double
precision numbers - so that the working interval contains of the
order of $10^{16}$ representable points - out of 1,000 randomly
chosen initial points (see Table \ref{tab:tab1}),
\begin{itemize}
  \item 596, i.e., the majority, converged to the fixed point
corresponding to the unstable fixed point $\{0\}$ in equation
\ref{eq2},
  \item 404 converged to a cycle of period 15,784,521.
\end{itemize}

Thus, in this case at least, the very long-term behaviour of
numerical orbits is, for a substantial fraction of initial points,
in flagrant disagreement with the true behaviour of typical orbits
of the original smooth logistic map.

In others numerical experiments we have performed, the computer
working with fixed finite precision is able to represent finitely
many points in the interval in question. It is probably good, for
purposes of orientation, to think of the case where the
representable points are uniformly spaced in the interval. The true
logistic map is then \emph{approximated} by a discretized map,
sending the finite set of representable points in the interval to
itself.

Describing the discretized mapping exactly is usually complicated,
but it is \emph{roughly} the mapping obtained by applying the exact
smooth mapping to each of the discrete representable points and
"rounding" the result to the nearest representable point. In our
experiments uniformly spaced points in the interval with several
order of discretization (ranging from 9 to 2,001 points) are
involved. In each experiment the questions addressed are:

\begin{itemize}
  \item how many periodic cycles are there and what are their periods ?
  \item how large are their respective basins of attraction, i.e. , for
each periodic cycle, how many initial points give orbits with
eventually land on the cycle in question ?
\end{itemize}

\begin{table}[!bth]
\begin{tabular}{rcr}
  \hline
  Period & Orbit & Relative Basin size \\
  \hline
  1 & $\{ 0\}$ (unstable fixed point) & 596 over 1,000 \\
  15,784,521 & Scattered over the interval & 404 over 1,000 \\
    \hline
\end{tabular}
\caption{Coexisting periodic orbits found using 1,000 random initial
points for double precision numbers} \label{tab:tab1}
\end{table}

\begin{table}[!bth]
\begin{tabular}{rrcr}
  \hline
  $N$ & Period & Orbit & Basin size \\
  \hline
  9 & 1 & $\{ 0\}$ & 3 over 9\\
  9 & 1 & $\{ 6\}$ & 2 over 9\\
  9 & 1 & $\{ 3,7\}$ & 4 over 9\\
   10 & 1 & $\{ 0\}$ & 2 over 10\\
   10 & 2 & $\{ 3,8\}$ & 8 over 10\\
 11 & 1 & $\{ 0\}$ & 3 over 11\\
 11 & 4& $\{ 3,8,6,9\}$ & 8 over 11\\
  \hline
\end{tabular}
\caption{Coexisting periodic orbits for the discretization with
regular meshes of $N=9$, $10$ and $11$ points} \label{tab:tab2}
\end{table}


\begin{table}[!bth]
\begin{tabular}{rrcr}
  \hline
  $N$ & Period & Orbit & Basin size \\
  \hline
  99 & 1 & $\{ 0\}$ & 3 over 99\\
  99 & 10 & $\{3,11,39,93,18,58,94,15,50,97\}$ & 96 over 99\\
  100 & 1 & $\{ 0\}$ & 2 over 100\\
   100 & 1 & $\{ 74\}$ & 2 over 100\\
   100 & 6 & $\{ 11,39,94,18,58,96\}$ & 72 over 100\\
 100  & 7 & $\{ 7,26,76,70,82,56,97\}$ & 24 over 100\\
 101 & 1& $\{0\}$ & 3 over 101\\
 101 & 1& $\{75\}$ & 2 over 101\\
 101 & 1& $\{16,61,95\}$ & 96 over 101\\
  \hline
\end{tabular}
\caption{Coexisting periodic orbits for the discretization with
regular meshes of $N=99$, $100$ and $101$ points} \label{tab:tab3}
\end{table}

\begin{table}[!bth]
\begin{tabular}{rrcr}
  \hline
  $N$ & Period & Orbit & Basin size \\
  \hline
  1,999 & 1 & $\{ 0\}$ & 3 over 1,999\\
  1,999 & 4 & $\{554;1,601;1,272;1,848\}$ &990 over 1,999\\
  1,999  & 8 & $\{3;11;43; 168; 615; 1,702;1,008;1,997\}$ & 1,006 over 1,999\\
   2,000 & 1 & $\{ 0\}$ & 2 over 2,000\\
   2,000 & 1 & $\{ 1,499\}$ & 14 over 2,000\\
 2,000  & 2& $\{691;1,808\}$ & 138 over 2,000\\
  2,000 & 3& $\{276;1,221;1,900\}$ & 6 over 2,000\\
  2,000 & 8& $\{3;11;43;168;615;1,703;1,008;1,998\}$ & 1,840 over 2,000\\
 2,001 & 1 & $\{ 0\}$ & 5 over 2,001\\
   2,001 & 1 & $\{1,500\}$ &34 over 2,001\\
   2,001  & 2 & $\{691;1,809\}$ & 92 over 2,001\\
    2,001 & 8 & $\{ 3;11;43;168;615;1,703;1,011;1,999\}$ & 608 over 2,001\\
    2,001& 18 & $\{35;137;510;1,519;1,461; 1,574; \ldots\}$ & 263 over 2,001\\
  2,001  & 25& $\{27;106; 401; 1,282; 1,840; 588; \ldots\}$ & 1,262 over 2,001\\
    \hline
\end{tabular}
\caption{Coexisting periodic orbits for the discretization with
regular meshes of $N=1,999$;$2,000$ and $2,001$ points}
\label{tab:tab4}
\end{table}

On an another hand, for relatively coarse discretizations the orbit
structure is determined completely, i.e., all the periodic cycles
and the exact sizes of their basins of attraction are found. Some
representative results are given in Tables \ref{tab:tab2} to
\ref{tab:tab4}. In theses tables, $N$ denotes the order of the
discretization, i.e., the representable points are the numbers,
$\frac{j}{N}$, with $0 \leq j < N$.

 The Table \ref{tab:tab2} shows coexisting periodic orbits for the
discretization with regular meshes of $N=9$, $10$ and $11$ points.
There are exactly 3, 2 and 2 cycles.\\
The Table \ref{tab:tab3} shows coexisting periodic orbits for the
discretization with regular meshes of $N=99$, $100$ and $101$
points. There are
exactly 2, 4 and 3 cycles.\\
The Table \ref{tab:tab4} shows coexisting periodic orbits for the
discretization with regular meshes of $N=1,999$, $N=2,000$ and
$N=2,001$ points.

It seems that the computation of numerical approximations of the
periodic orbits leads to unpredictable results.

\subsection{Statistical properties}

Many others examples could be given, but those given may serve to
illustrate the intriguing character of the results: the outcomes
proves to be extremely sensitive to the details of the experiment,
but the results all have a similar flavour : a relatively small
number of cycles attract near all orbits, and the lengths of these
significant cycles are much larger than one but much smaller than
the number of representable points.

In \cite{Dia}, P. Diamond and A. Pokrovskii, suggest that
statistical properties of the phenomenon of computational collapse
of discretized chaotic mapping can be modelled by random mappings
with an absorbing centre. The model gives results which are very
much in line with computational experiments and there appears to be
a type of universality summarised by an Arcsine law. The effects are
discussed with special reference to the family of mappings

\begin{equation}\label{eq4}
x_{n+1}=1-|1-2x_n|^{\ell}~~~~0\leq x \leq 1 ~~~ 1 \leq \ell \leq 2
\end{equation}

Computer experiments show close agreement with prediction of the
model.

However these results are of statistical nature, they do not give
accurate information on the exact nature of the orbits (e.g. length
of the shortest one, of the greater one, size of their basin of
attraction ...). It is why we consider the problem of computational
discrepancies in an original way in the next section.

\section{Pattern defined by all the orbits of a dynamical system}\label{section3}

In this section in order to describe more precisely which kind of
behaviour occurs in discretized dynamical systems on finite sets we
conceive a new mathematical tool: the global orbit pattern of a
function that is the set of the periods of every different orbits of
the dynamical system associated to the function when the initial
points are took in increasing order.

\subsection{General definitions}
For every $x_0 \in X$, let $\lbrace x_i \rbrace$ be the sequence of
the orbit of the dynamical system associated to the function $f$
which maps $X$ onto $X$ defined by

\begin{equation}\label{eq5}
    x_{i+1}=f(x_i) ~~\textrm{for}~~~i \geq
0 .
\end{equation}

For convenience $\forall x_0 \in X$ we denote

\begin{equation}\label{eq6}
  f^0(x_0)=x_0
\end{equation}

and
\begin{equation}\label{eq7}
 \forall p \geq 1,\forall x_0 \in X,~~f^p(x_0)=\underbrace{f \circ f
 \circ \ldots \circ f }_{p ~~\textrm{times}}(x_0).
\end{equation}

 Hence
\begin{equation}\label{eq8}
 x_i=f^i(x_0).
\end{equation}

  The orbit of $x_0$ under $f$ is the set of points
  $\mathcal{O}(x_0,f)=\lbrace  f^i(x_0), i \geq 0 \rbrace = \lbrace x_i, i \geq 0 \rbrace$.\\
The starting point $x_0$ for the orbit is called the initial value
of the orbit.\\
A point $x$ is a fixed point of the map $f$ if $f(x)=x$.\\
A point $x$ is a periodic point with period $p$
if $f^p(x)=x$ and $f^k(x)\neq x$ for all $k$ such that $0 \leq k <p$, $p$ is called the order of $x$.\\
If $x$ is periodic of order $p$, then the orbit of $x$ under $f$ is
the finite set \\$\lbrace x, f(x), f^2(x), \ldots, f^{p-1}(x)
\rbrace$. We will call this set the periodic orbit of order $p$
or a $p$-cycle.\\
A fixed point is then a 1-cycle.\\
The point $x$ is an eventually periodic point of $f$ with order $p$
if there exists $K>0$ such that $\forall k\geq K$
$f^{k+p}(x)=f^k(x)$.\\
$\forall x \in X$, we denote $\omega(x,f)$ the order of $x$
under $f$ or simply $\omega(x)$ when the map $f$ involved is obvious.\\
A subset $T$ of $X$ is invariant under $f$ if $f^{-1}(T)=T$. That is
equivalent to say that $T$ is invariant under $f$ if and only if
$f(T) \subset T$ and $f^{-1}(T)\subset T$.\\

\nota{$\sharp X$ is the cardinal of the finite set $X$.}

\subsection{Map on finite set}
Along this paper, $N$ is a non-zero integer and $\sharp A$ stands
for the cardinal of any finite set $A$. In this article, we consider
$X$ as an ordered finite set with $N$ elements. We denote it $X_N$,
it is isomorphic to the interval $\dbl 0, N-1 \dbr \subset
\mathbb{N}$. Then $\sharp X_N=N$. Let $f$ be a map from $X_N$ into
$X_N$. We denote by $\mathcal{F}_N$ the set of the maps from $X_N$
into $X_N$. Clearly, $\mathcal{F}_N$ is a finite set and $\sharp
\mathcal{F}_N= N^N$ elements. For all $x \in X_N$,
$\mathcal{O}(x,f)$ is necessarily a finite set with at most $N$
elements. Indeed, let us consider the sequence $\lbrace x, f(x),
f^2(x), \ldots , f^{N-1}(x), f^N(x)\rbrace$ of the first $N+1$
iterated points. Thanks to the Dirichlet's box principle, two
elements are equals because $X_N$ has exactly $N$ different values.
Thus, every initial value of $X_N$ leads ultimately to a repeating
cycle. More precisely, if $x$ is a fixed point $\mathcal{O}(x,f)$ is
the unique element $x$ and if $x$ is a periodic point with order
$p$, $\mathcal{O}(x,f)$ has exactly $p$ elements. In this case, the
orbit of $x$ under $f$ is the set $\mathcal{O}(x,f)= \lbrace x,
f(x), f^2 (x), \ldots, f^{p-1}(x) \rbrace$. If $x$ is an eventually
periodic point with order $p$, there exists $K>0$ such that $\forall
k\geq K$ $f^{k+p}(x)=f^k(x)$. In this case, the orbit of $x$ under
$f$ is the set $\mathcal{O}(x,f)= \lbrace x, f(x), f^2(x), \ldots,
f^K(x), f^{K+1}(x), \ldots, f^{K+p-1}(x) \rbrace $.

\subsection{Equivalence classes}

 \subsubsection{Components}
Let $f \in \mathcal{F}_N$. We consider on $X_N$  the relation $\sim$
defined by : $\forall x, x' \in X_N$,~ $x \sim x' \Leftrightarrow
\exists k \in \mathbb{N}$ such that $f^k(x) \in \mathcal{O}(x',f)$.
The relation $\sim$ is an equivalence relation on $X_N$.
$\mathcal{S}_N/\sim$ is the collection of the equivalence classes
that we will call components of $X_N$ under $f$ which constitute a
partition of $X_N$. The number of components are given in
\cite{Krus}. Asymptotic properties of the number of cycles and
components are studied in \cite{Mutaf}. For each component, we take
as class representative element the least element of the component.
The components will be written $T_N(x_0,f), \ldots,
T_N(x_{p_{f,N}},f)$ where
$x_i$ is the least element of $T_N(x_i,f)$. \\

By analogy with real dynamical systems, we can define attractive
and repulsive components in discretized dynamical systems as follows.\\

\defn{A component is repulsive when it is a cycle. Otherwise, the
component is attractive.}

\rem{In other words, a component is attractive when the component
contains at least an eventually periodic element. The corresponding
cycle is strictly contained in an attractive component.}

Examples are given in Tables \ref{tab:tab5}, \ref{tab:tab6} and
\ref{tab:tab7}.\\

For instance, in Table \ref{tab:tab6}, the fonction $f$ has
$\{2,7\}$ as period-2 orbit and $\{1,2,7,9\}$ as component which is
attractive because $1$ and $2$ are eventually periodic elements. \\

\begin{table}[!bth]
\begin{tabular}{lc}
  \hline
  Function & orbit/component/nature  \\
    \hline

 \begin{tabular}{ccc}
 0 & $\rightarrow$ & 6\\
  1 & $\rightarrow$ & 3\\
   2 & $\rightarrow$ & 2\\
    3 & $\rightarrow$ & 5\\
     4 & $\rightarrow$ & 8\\
      5 & $\rightarrow$ & 10\\
      6 & $\rightarrow$ & 9\\
       7 & $\rightarrow$ & 4\\
        8 & $\rightarrow$ & 7\\
         9 & $\rightarrow$ & 6\\
  10 & $\rightarrow$ & 5\\
\end{tabular}
 &
\begin{tabular}{ccc}
 period-2 orbit : $\{6,9\}$ &$\{0,6,9\}$ & attractive\\
  & & \\
 period-2 orbit : $\{5,10\}$ & $\{1,3,5,10\}$ & attractive\\
  & & \\
 fixed point : $\{2\}$ & $\{2\}$ & repulsive\\
   & & \\
  period-3 orbit : $\{4,8,7\}$ & $\{4,8,7\}$ & repulsive \\
\end{tabular}\\
\hline

\end{tabular}
\caption{Orbits and components of a function belonging to
$\mathcal{F}_{11}$ with gop $[2,2,1,3]_{11}$.} \label{tab:tab5}
\end{table}


\begin{table}[!bth]
\begin{tabular}{lc}
  \hline
  Function & orbit/component/nature  \\
    \hline

 \begin{tabular}{ccc}
 0 & $\rightarrow$ & 4\\
  1 & $\rightarrow$ & 2\\
   2 & $\rightarrow$ & 7\\
    3 & $\rightarrow$ & 3\\
     4 & $\rightarrow$ & 8\\
      5 & $\rightarrow$ & 10\\
      6 & $\rightarrow$ & 5\\
       7 & $\rightarrow$ & 2\\
        8 & $\rightarrow$ & 4\\
         9 & $\rightarrow$ & 1\\
          10 & $\rightarrow$ & 6\\
\end{tabular}
 &
\begin{tabular}{ccc}
 period-2 orbit : $\{4,8\}$ &$\{0,4,8\}$ & attractive\\
   & & \\
 period-2 orbit : $\{2,7\}$ & $\{1,2,7,9\}$ & attractive\\
  & & \\
 fixed point : $\{3\}$ & $\{3\}$ & repulsive\\
   & & \\
  period-3 orbit : $\{5,10,6\}$ & $\{5,10,6\}$ & repulsive \\
\end{tabular}\\
\hline

\end{tabular}
\caption{Orbits and components of a function belonging to
$\mathcal{F}_{11}$ with gop $[2,2,1,3]_{11}$.} \label{tab:tab6}
\end{table}

\begin{table}[!bth]
\begin{tabular}{lc}
  \hline
  Function & orbit/component/nature  \\
    \hline

 \begin{tabular}{ccc}
 0 & $\rightarrow$ & 9\\
  1 & $\rightarrow$ & 6\\
   2 & $\rightarrow$ & 4\\
    3 & $\rightarrow$ & 7\\
     4 & $\rightarrow$ & 10\\
      5 & $\rightarrow$ & 3\\
      6 & $\rightarrow$ & 1\\
       7 & $\rightarrow$ & 5\\
        8 & $\rightarrow$ & 2\\
         9 & $\rightarrow$ & 0\\
            10 & $\rightarrow$ & 10\\
\end{tabular}
 &
\begin{tabular}{ccc}
  period-2 orbit : $\{0,9\}$ &$\{0,9\}$ & repulsive\\
   & & \\
 period-2 orbit : $\{1,6\}$ & $\{1,6\}$ & repulsive\\
  & & \\
 fixed point : $\{10\}$ & $\{2,4,8,10\}$ & attractive\\
   & & \\
  period-3 orbit : $\{3,7,5\}$ & $\{3,7,5\}$ & repulsive \\
\end{tabular}\\
\hline

\end{tabular}
\caption{Orbits and components of a function belonging to
$\mathcal{F}_{11}$ with gop $[2,2,1,3]_{11}$.} \label{tab:tab7}
\end{table}

\subsubsection{Order of elements}

Here are some remarks on the order of elements of components.\\

\rem{The order of every element of a component is the length of its
inner cycle.}

\defn{For all $x \in X_N$, there exists $i \in \dbl 0, p_{f,N} \dbr $ such that $x$ belongs to the component $T_N(x_i,f)$.
Then $\omega(x,f)$ is equal to the order $\omega(x_i,f)$. }

\rem{For all $i \in \dbl 0,p_{f,N} \dbr$, $T_N(x_i,f)$ is an
invariant subset of $X_N$ under $f$.}

In the example given in Table \ref{tab:tab5}, the order of the
element $0$ is 2, the order of the element $1$ is $2$, the order of
the element $4$ is 3. The elements $1$ and $3$ have the same order.

\subsection{Definition of global orbit pattern}
For each $f \in \mathcal{F}_N$, we can determine the components of
$X_N$ under $f$. For each component, we determine the order of any
element. Thus, for each $f \in \mathcal{F}_N$, we have a set of
orders that we will denote $\Omega(f, N)$. Be given $f$, there exist
$p_{f,N}$ components and $p_{f,N}$ representative elements such that
$x_0 < x_1 < \ldots <x_{p_{f,N}}$.\\

For each $f \in \mathcal{F}_N$, the sequence $[\omega(x_0),
\omega(x_1), \ldots, \omega(x_{p_{f,N}}); f]_{\mathcal{F}_N}$ with
$x_0 < x_1 < \ldots <x_{p_{f,N}}$ will design
 the global orbit pattern of $f \in \mathcal{F}_N$. \\

We will write $gop(f)=[\omega(x_0),
\omega(x_1), \ldots, \omega(x_{p_{f,N}}); f]_{\mathcal{F}_N}$.\\

When the reference to $f \in \mathcal{F}_N$ is obvious, we will
write shortly \\$gop(f)=[\omega(x_0), \omega(x_1), \ldots,
 \omega(x_{p_{f,N}})]_N$ or $gop(f)=[\omega(x_0),
\omega(x_1), \ldots, \omega(x_{p_{f,N}})]$ .\\

For example, the same gop associated to the functions given in
Tables \ref{tab:tab5}, \ref{tab:tab6} and
\ref{tab:tab7} is $[2,2,1,3]_{11}$.\\

Another example is given in Table \ref{tab:tab8}.
In that example, we have $\omega(0)=2$, $\omega(3)=1$, $\omega(4)=4$. \\\\

\begin{table}[!bth]
\begin{tabular}{lc}
  \hline
  Function & orbit/component/nature  \\
    \hline

 \begin{tabular}{ccc}
 0 & $\rightarrow$ & 1\\
  1 & $\rightarrow$ & 0\\
   2 & $\rightarrow$ & 0\\
    3 & $\rightarrow$ & 3\\
     4 & $\rightarrow$ & 5\\
      5 & $\rightarrow$ & 6\\
      6 & $\rightarrow$ & 7\\
       7 & $\rightarrow$ & 4\\
\end{tabular}
 &
\begin{tabular}{ccc}
period-2 orbit : $\{0,1\}$ &$\{0,1,2\}$ & attractive\\
   & & \\
fixed point : $\{3\}$ & $\{3\}$ & repulsive\\
  & & \\
 period-4 orbit : $\{4,5,6,7\}$ & $\{4,5,6,7\}$ & repulsive\\
& & \\

\end{tabular}\\
\hline

\end{tabular}
\caption{Orbits and components of a function belonging to
$\mathcal{F}_{8}$ with gop $[2,1,4]_{8}$.} \label{tab:tab8}
\end{table}

\defn{The set of all the global orbit patterns of $\mathcal{F}_N$ is called $\mathcal{G}(\mathcal{F}_N)$.}

For example, for $N=5$, the set $\mathcal{G}(\mathcal{F}_5)$ is\\
\begin{minipage}{14cm}
$\{ [1]; [1,1] ;[1,1,1] ; [1,2] ; [1,1,1,1] ; [1,1,2]\linebreak
[1,2,1]; [1,3]; [1,1,1,1,1] ; [1,1,1,2] ; [1,1,2,1] ;\linebreak
[1,1,3] ; [1,2,1,1]; [1,2,2] ; [1,3,1] ; [1,4] ;\linebreak   [2] ;
[2,1] ; [2,1,1]; [2,2] ; [2,1,1,1] ; [2,1,2] ; [2,2,1]\linebreak
[2,3] ; [3]; [3,1]; [3,1,1] ; [3,2]; [4] ; [4,1];[5]\}$.
\end{minipage}

\subsection{Class of gop}
We give the following definitions : \\

\defn{Let be $A=[\omega_1, \ldots ,\omega_p]_N$ a gop. Then the class of $A$, written $\overline{A}$, is the set of all the
functions $f \in \mathcal{F}_N$ such that the global orbit pattern
associated to $f$ is $A$.}

For example, for $N=11$, the class of the gop
$\overline{[2,2,1,3]}_{11}$ contains the following few of many
functions defined in Tables \ref{tab:tab5}, \ref{tab:tab6} and
\ref{tab:tab7}. The periodic orbit which are encountered have the same length nevertheless there are different.\\

\defn{Let be $A=[\omega_1, \ldots ,\omega_p]_N$ a gop.
\\Then the modulus of $A$ is $|A| = \sum\limits_{k=1}^p \omega_k$.}

\rem{$\left|[\omega_1, \ldots ,\omega_p]_N\right| \leq N$.}

\nota{ $[\omega_{\widetilde{k}}]_N$ means $[\underbrace{\omega,
\ldots, \omega}_{k ~\textrm{times}}]_N$ and
$[\omega_{\widetilde{k}},\nu_{\widetilde{m}}]_N$ means
$[\underbrace{\omega, \ldots, \omega}_{k
~\textrm{times}},\underbrace{\nu, \ldots, \nu}_{m ~\textrm{times}}
]_N$.}

\subsection{Threshold functions}

\subsubsection{Ordering the discrete maps}

\thm{The sets $\mathcal{F}_N$ and $\dbl 1,N^N \dbr$ are isomorphic.}

\proof{ We define the function $\phi$ from $\mathcal{F}_N$ to $\dbl
1,N^N \dbr$ by : for each  $f \in \mathcal{F}_N$, $\phi(f)$ is the
integer $n$ such that $n = \sum\limits_{k=0}^{N-1} f(k)N^{N-1-k} + 1$.\\
Then
$\phi$ is well defined because $n \in \dbl 1,N^N \dbr$. \\
Let $n$ be a given integer between 1 and $N^N$. We convert $n-1$ in
base $N$ : there exists a unique $N$-tuple $(a_{n-1,0}; a_{n-1,1};
\ldots ; a_{n-1,N-1}) \in \dbl 0,N-1\dbr^N$ such that
$\overline{n-1}^N = \sum\limits_{i=0}^{N-1} a_{n-1,N-1-i}
N^{N-i-1}$. We can thus define the map $f_n$ with : $\forall i \in
X_N$, $f_n(i)=a_{n-1,N-i-1}$. Then $\phi$ is one to one. }\\

\rem{This implies $\mathcal{F}_N$ is totally ordered.}

\defn{Let $f \in \mathcal{F}_N$.
Then
\begin{equation}\label{eq9}
n=\sum\limits_{k=0}^{N-1} f(k)N^{N-1-k}+ 1
\end{equation}
is called
the rank of $f$.  }

\subsubsection{Threshold functions} Be given a global orbit pattern $A$, we
are exploring the class $\overline{A}$.\\

\thm{For every $A \in \mathcal{G}(\mathcal{F}_N)$, the class
$\overline{A}$ has a unique function with minimal rank.}

\defn{For every class $\overline{A} \in \mathcal{G}(\mathcal{F}_N)$, the function
defined by the previous theorem will be called the threshold
function for the class $\overline{A}$ and will be denoted by $Tr(A)$
or $Tr(\overline{A})$. }

To prove the theorem, we need the following definition : \\

\defn{Let $f \in \mathcal{F}_N$ be a function. Let $p$ a non zero integer smaller than $N$. Let be $x_1, \ldots, x_{p}$
$p$ consecutive elements of $X_N$. Then $x_1, \ldots, x_{p}$ is a
canonical $p$-cycle in relation to $f$ if  $\forall j \in \dbl 1,
p-1 \dbr$ , $f(x_j)=x_{j+1}$ and $f(x_p)=x_1$.}

\proof{Let $[\omega_1, \ldots ,\omega_p]$ be a global orbit pattern
of $\mathcal{G}(\mathcal{F}_N)$. We construct a specific function
$f$ belonging to the class $\overline{[\omega_1, \ldots ,\omega_p]}$
and we prove that the function so obtained is the smallest with
respect to the order on $\mathcal{F}_N$. With the first $\omega_1$
elements of $\dbl 0,N-1 \dbr$, that is the set of integers $\dbl
0,\omega_1-1 \dbr$, we construct the canonical $\omega_1$-cycle : if
$\omega_1=1$, we define $f(0)=0$, else $f(0)=1$, $f(1)=2$, $\ldots$,
$f(\omega_1-2)=\omega_1-1$, $f(\omega_1-1)=0$.
\\Then $\forall j \in \dbl \omega_1-1, \omega_1+N-s-1 \dbr$, we define
$f(j)=0$. \\Then with the next $\omega_2$ integers $\dbl
\omega_1+N-s,\omega_1+N-s+\omega_2-1 \dbr$ we construct the
canonical $\omega_2$-cycle. We keep going on constructing for all
 $ j \in  \dbl 3,p \dbr$ the canonical $\omega_j$-cycle.\\
 In consequence, we have found a function $f$ belonging to the class $\overline{[\omega_1, \ldots ,\omega_p]}$. \\
Assume there exists a function $g \in \mathcal{F}_N$ belonging to
the class of $f$ such that $g < f$. Let $I= \{i \in \dbl 0,N-1 \dbr
~\textrm{such that}~ f(i) \neq 0 \} $. As $g <f$, there exists $i_0
\in I$ such that $g(i_0) <f(i_0)$. There exists also $j_0$ such that
$i_0 \in \omega_{j_0}$. If $f(i_0)=i_0$, then $\omega_{j_0}=1$,
$g(i_0)<i_0$ and then $g(i_0) \notin \omega_{j_0}$. Then the global
orbit pattern of $g$ doesn't contain anymore 1 as cycle. The global
orbit pattern of $g$ is different from the global orbit pattern of
$f$. If $f(i_0)=i_0+1$, then $g(i_0)\leq i_0$. Either $g(i_0)=i_0$
and then the global orbit pattern of $g$ is changed, or $g(i_0)<i_0$
and we are in the same situation as previously. Thus, in any case,
the smallest function belonging to the class $[\overline{\omega_1,
\ldots ,\omega_p}]$ is the
one constructed in the first part of the proof.}\\

The proof of the theorem gives an algorithm of construction of the
threshold function associated to a given gop.\\
The threshold function associated to the gop
$[2_{\widetilde{2}},1,3]_{11}$ is
explained in Table \ref{tab:tab9}. Its rank is $n=25,938,474,637$.\\

\begin{table}[!bth]
\begin{tabular}{ccccc}
  \hline
  First step & Second step & Third step & Fourth step & Fifth step  \\
    \hline
    \begin{minipage}{2,4cm}
Construction of the first canonical 2-cycle\\
 \begin{tabular}{ccc}
 0 & $\rightarrow$ & 1\\
  1 & $\rightarrow$ & 0\\
   2 & $\rightarrow$ & \\
    3 & $\rightarrow$ & \\
     4 & $\rightarrow$ & \\
      5 & $\rightarrow$ & \\
      6 & $\rightarrow$ & \\
       7 & $\rightarrow$ & \\
        8 & $\rightarrow$ & \\
         9 & $\rightarrow$ & \\
           10 & $\rightarrow$ & \\
\end{tabular}
\end{minipage}
 &
     \begin{minipage}{2,4cm}
Construction of the last canonical 3-cycle\\
 \begin{tabular}{ccc}
 0 & $\rightarrow$ & 1\\
  1 & $\rightarrow$ & 0\\
   2 & $\rightarrow$ & \\
    3 & $\rightarrow$ & \\
     4 & $\rightarrow$ & \\
      5 & $\rightarrow$ & \\
      6 & $\rightarrow$ & \\
       7 & $\rightarrow$ & \\
        8 & $\rightarrow$ & 9\\
         9 & $\rightarrow$ & 10\\
            10 & $\rightarrow$ & 8\\
\end{tabular}
\end{minipage}
&

    \begin{minipage}{2,4cm}
   Construction of the canonical 1-cycle\\
 \begin{tabular}{ccc}
 0 & $\rightarrow$ & 1\\
  1 & $\rightarrow$ & 0\\
   2 & $\rightarrow$ & \\
    3 & $\rightarrow$ & \\
     4 & $\rightarrow$ & \\
      5 & $\rightarrow$ & \\
      6 & $\rightarrow$ & \\
       7 & $\rightarrow$ & 7\\
         8 & $\rightarrow$ & 9\\
         9 & $\rightarrow$ & 10\\
            10 & $\rightarrow$ & 8\\
\end{tabular}
\end{minipage}
&

    \begin{minipage}{2,4cm}
Construction of the canonical 2-cycle\\
 \begin{tabular}{ccc}
 0 & $\rightarrow$ & 1\\
  1 & $\rightarrow$ & 0\\
   2 & $\rightarrow$ & \\
    3 & $\rightarrow$ & \\
     4 & $\rightarrow$ & \\
      5 & $\rightarrow$ & 6\\
      6 & $\rightarrow$ & 5\\
       7 & $\rightarrow$ & 7\\
         8 & $\rightarrow$ & 9\\
         9 & $\rightarrow$ & 10\\
            10 & $\rightarrow$ & 8\\
\end{tabular}
\end{minipage}
&

\begin{minipage}{2,4cm}
Filling the remaining images with 0\\
 \begin{tabular}{ccc}
0 & $\rightarrow$ & 1\\
  1 & $\rightarrow$ & 0\\
   2 & $\rightarrow$ & 0\\
    3 & $\rightarrow$ & 0\\
     4 & $\rightarrow$ & 0\\
      5 & $\rightarrow$ & 6\\
      6 & $\rightarrow$ & 5\\
       7 & $\rightarrow$ & 7\\
         8 & $\rightarrow$ & 9\\
         9 & $\rightarrow$ & 10\\
            10 & $\rightarrow$ & 8\\
\end{tabular}
\end{minipage}\\
\hline

\end{tabular}
\caption{Algorithm for the threshold function construction for the
gop $[2_{\widetilde{2}},1,3]_{11}$.} \label{tab:tab9}
\end{table}

\thm{There are exactly $2^N-1$ different global orbit patterns in
$\mathcal{F}_N$.}

That is

\begin{equation}\label{eq10}
\sharp \mathcal{G}(\mathcal{F}_N)=2^N-1.
\end{equation}\\

For example, for $N=4$, $\sharp \mathcal{G}(\mathcal{F}_4)=2^4-1=15$.\\

\proof{ Let $p$ an integer between 1 and $N$. Consider the set
$L(p,N)$ of $p$-tuples $(\alpha_1, \ldots, \alpha_p) \in
(\mathbb{N}^*)^p $ such that $\alpha_1+ \ldots+ \alpha_p \leq N$.
\\We write $L(N) = \lbrace L(p,N), p=1 \ldots N \rbrace$. $L(N)$ and
$\mathcal{G}(\mathcal{F}_N)$ have the same elements. Then \\$
\displaystyle \sharp \mathcal{G}(\mathcal{F}_N)= \sum
\limits_{p=1}^{p=N} \sharp L(p,N)= \sum \limits_{p=1}^{p=N}
\left(\begin{array}{c}
  N \\
  p \\
\end{array}\right)=2^N-1$.
}

\subsubsection{Ordering the global orbit patterns}

We define an order relation on $\mathcal{G}(\mathcal{F}_N)$.\\

\prop{Let $A$ and $B$ be two global orbit patterns of
$\mathcal{G}(\mathcal{F}_N)$.\\ We define the relation $\prec$ on
the set $\mathcal{G}(\mathcal{F}_N)$ by $$A \prec B~~ \textit{iff}
~~Tr(A)<Tr(B)$$ Then the set $(\mathcal{G}(\mathcal{F}_N),\prec)$ is
totally ordered.}

\proof{As the order  $\prec$ refers to the natural order of
$\mathbb{N}$, the proof is obvious.}\\

Let $r \geq 1$, $p\geq 1$ be two integers. Let $[\omega_1, \ldots
,\omega_p]$ and $[\omega'_1, \ldots, \omega'_r]$ be two global orbit
patterns of $\mathcal{G}(\mathcal{F}_N)$. For example, if $p<r$, in
order to compare them, we admit that we can fill $[\omega_1, \ldots
,\omega_p]$ with $r-p$ zeros and write $[\omega_1, \ldots
,\omega_p]=[\omega_1, \ldots ,\omega_p,
0, \ldots, 0]$.\\

\prop{Let $r \geq 1$, $p\geq 1$ be two integers such that $p \leq
r$. Let $A =[\omega_1, \ldots ,\omega_p]$ and $B =[\omega'_1,
\ldots, \omega'_r]$ be two global orbit patterns. \\

\begin{minipage}{14cm}
\begin{itemize}
\item If $r=p=1$ and $\omega_1 < \omega'_1$ then $A \prec B$.
\item If $r \geq 2$ then
\begin{itemize}
  \item[$\ast$]  If $\omega_1 < \omega'_1$ then $A \prec B$.
    \item[$\ast$]  If $\omega_1 = \omega'_1$ then there exists $K \in \dbl 2;r \dbr$
such that $\omega_K \neq \omega'_K$ and $\forall i < K$, $\omega_i =
\omega'_i$.
\begin{itemize}
  \item[$\bullet$] If $|A| < |B|$, then
$A \prec B$.
   \item[$\bullet$] If $|A| = |B|$, then if
$\omega_K < \omega'_K$ then $A \prec B$.
\end{itemize}
\end{itemize}
\end{itemize}
\end{minipage} \\}

For example, for $N=5$, the global orbit patterns are in increasing
order : $ [1]\prec [1_{\widetilde{2}}] \prec [1_{\widetilde{3}}]
\prec [1,2] \prec[1_{\widetilde{4}}] \prec [1_{\widetilde{2}},2]
\prec [1,2,1] \prec [1,3] \prec[1_{\widetilde{5}}]  \prec
[1_{\widetilde{3}},2] \prec [1_{\widetilde{2}},2,1] \prec
[1_{\widetilde{2}},3] \prec [1,2,1_{\widetilde{2}}] \prec
[1,2_{\widetilde{2}}] \prec[1,3,1] \prec[1,4] \prec [2] \prec [2,1]
\prec [2,1_{\widetilde{2}}] \prec[2_{\widetilde{2}}] \prec
[2,1_{\widetilde{3}}] \prec [2,1,2]
 \prec [2_{\widetilde{2}},1]  \prec [2,3]  \prec [3] \prec [3,1] \prec [3,1_{\widetilde{2}}]  \prec [3,2] \prec [4]  \prec
[4,1]\prec [5]$.

\subsubsection{Algorithm for ordering the global orbit patterns : a pseudo-decimal order}

The Table \ref{tab:tab10} gives a method for ordering the gop :
indeed, we consider each gop as if each one represents a decimal
number : we begin to order them in considering the first order
$\omega_1$. Considering two gops $A= [\omega_1, ..., \omega_p]$ and
$A'=[\omega_1', ..., \omega_r']$, if $\omega_1<\omega_1'$, then $A
\prec A'$. For example, $[2,1,2] \prec [4,1]$. If
$\omega_1=\omega_1'$ and $|A|-\omega_1<|A'|-\omega_1'$, then $A
\prec A'$. For example to compare the gop $[1,2]$ and the gop
$[1_{\widetilde{4}}]$, we say that the first order $\omega_1$ stands
for the unit digit - which is $\omega_1=1$ here, then the decimal
digits are respectively $0.2$ and $0.111$. We calculate for each of
them the modulus-$\omega_1$ : we find $|[1,2]| -1 = 2$ and
$|[1_{\widetilde{4}}]| -1 = 3$, thus $[1,2] \prec
[1_{\widetilde{4}}]$. Finally, if $\omega_1=\omega_1'$ and
$|A|-\omega_1=|A'|-\omega_1'$, then also we use the order of the
decimal part. For example, $[1_{\widetilde{5}}] \prec
[1,2,1_{\widetilde{2}}$ because $1.1111 < 1.211$. Applying this
process, we have the sequence of the
ordered gop for $N=4$ given in the previous paragraph.\\

\begin{table}[!bth]
\begin{tabular}{lcc||lcc}
  \hline
  Gop & Modulus & Modulus-$\omega_1$ &Gop & Modulus & Modulus-$\omega_1$\\
  \hline
  $[1] $& 1 & 0 & $[2]$& 2&0 \\
    $[1_{\widetilde{2}}]$ & 2 &1 & $[2,1]$ &3&1\\
    $ [1_{\widetilde{3}}]$ & 3 &2 & $[2,1_{\widetilde{2}}]$&4&2\\
      $ [1,2]$ & 3  &2 & $[2_{\widetilde{2}}]$&4&2\\
       $ [1_{\widetilde{4}}]$ & 4 &3 &$[2,1_{\widetilde{3}}]$&5&3\\
       $ [1_{\widetilde{2}},2]$ & 4 &3&$[2,1,2]$&5&3\\
        $ [1,2,1]$ & 4& 3&$[2_{\widetilde{2}},1]$&5&3\\
       $ [1,3]$ & 4 &3 &$[2,3]$&5&3\\
       $[1_{\widetilde{5}}]$ & 5&4 &&\\
       $[1_{\widetilde{3}},2]$ &5&4&[3] &3&0\\
       $[1_{\widetilde{2}},2,1]$ &5&4&$[3,1]$&4&1\\
       $[1_{\widetilde{2}},3]$ &5&4&$[3,1_{\widetilde{2}}]$&5&2\\
       $[1,2,1_{\widetilde{2}}]$ &5&4&$[3,2]$&5&2\\
       $[1,2_{\widetilde{2}}]$ & 5&4&&\\
       $[1,3,1]$& 5&4&$[4]$&4&0\\
       $[1,4]$ & 5& 4&$[4,1]$&5&1\\
       &&&&&\\
       & & &$[5]$&5&1\\

  \hline
\end{tabular}
\caption{Ordered gop for $N=5$ with modulus and modulus-$\omega_1$}
\label{tab:tab10}
\end{table}

\begin{figure}
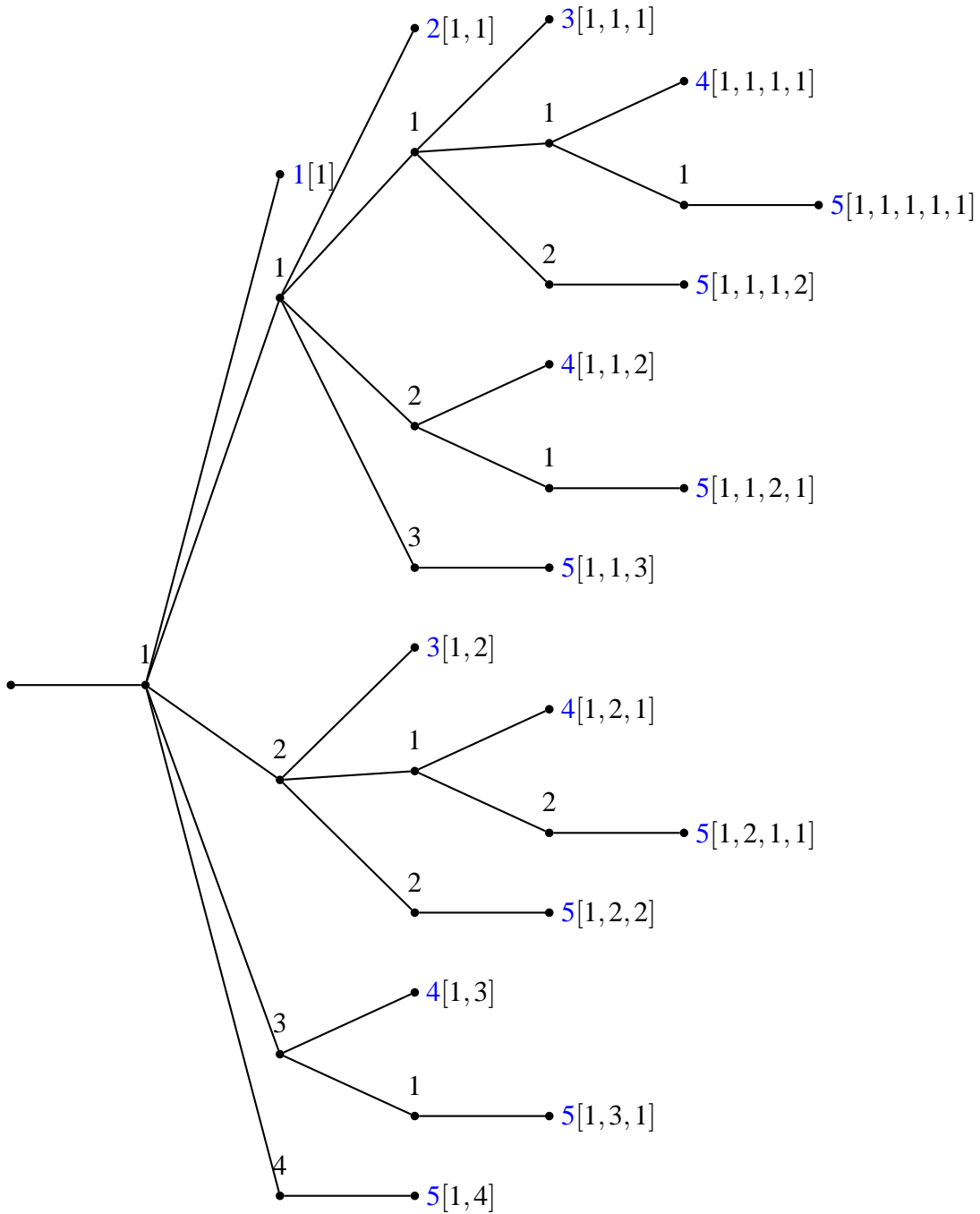

\psset{nodesep=0mm,levelsep=20mm,treesep=10mm}
\pstree[treemode=R]{\Tdot} { \pstree
{\Tdot~[tnpos=a]{$1$}\taput{\small $$}} {
\Tdot~[tnpos=r]{$\textcolor[rgb]{0.00,0.00,1.00}{1}[1]$}\taput{\small
$$} \pstree {\Tdot~[tnpos=a]{$1$}\taput{\small $$}} {
\Tdot~[tnpos=r]{$\textcolor[rgb]{0.00,0.00,1.00}{2}
[1,1]$}\taput{\small $$} \pstree {\Tdot~[tnpos=a]{$1$}\taput{\small
$$}} { \Tdot~[tnpos=r]{$\textcolor[rgb]{0.00,0.00,1.00}{3} [1,1,1]$}\taput{\small $$} \pstree
{\Tdot~[tnpos=a]{$1$}\taput{\small $$}} {
\Tdot~[tnpos=r]{$\textcolor[rgb]{0.00,0.00,1.00}{4}
[1,1,1,1]$}\taput{\small $$} \pstree
{\Tdot~[tnpos=a]{$1$}\tbput{\small $$}} {
\Tdot~[tnpos=r]{$\textcolor[rgb]{0.00,0.00,1.00}{5}
[1,1,1,1,1]$}\taput{\small $$} } } \pstree
{\Tdot~[tnpos=a]{$2$}\tbput{\small $$}} { \Tdot~[tnpos=r]{$
\textcolor[rgb]{0.00,0.00,1.00}{5} [1,1,1,2]$}\taput{\small $$} } }
\pstree {\Tdot~[tnpos=a]{$2$}\taput{\small $$}} {
\Tdot~[tnpos=r]{$\textcolor[rgb]{0.00,0.00,1.00}{4}
[1,1,2]$}\taput{\small $$} \pstree
{\Tdot~[tnpos=a]{$1$}\tbput{\small $$}} {
\Tdot~[tnpos=r]{$\textcolor[rgb]{0.00,0.00,1.00}{ 5}
[1,1,2,1]$}\taput{\small $$} } } \pstree
{\Tdot~[tnpos=a]{$3$}\taput{\small $$}} {
\Tdot~[tnpos=r]{$\textcolor[rgb]{0.00,0.00,1.00}{5}
[1,1,3]$}\taput{\small $$} } } \pstree
{\Tdot~[tnpos=a]{$2$}\taput{\small $$}} {
\Tdot~[tnpos=r]{$\textcolor[rgb]{0.00,0.00,1.00}{3}
[1,2]$}\taput{\small $$} \pstree {\Tdot~[tnpos=a]{$1$}\taput{\small
$$}} { \Tdot~[tnpos=r]{$\textcolor[rgb]{0.00,0.00,1.00}{4} [1,2,1]$}\taput{\small $$} \pstree
{\Tdot~[tnpos=a]{$2$}\tbput{\small $$}} {
\Tdot~[tnpos=r]{$\textcolor[rgb]{0.00,0.00,1.00}{5}
[1,2,1,1]$}\taput{\small $$} } } \pstree
{\Tdot~[tnpos=a]{$2$}\tbput{\small $$}} {
\Tdot~[tnpos=r]{$\textcolor[rgb]{0.00,0.00,1.00}{5}
[1,2,2]$}\taput{\small $$} } } \pstree
{\Tdot~[tnpos=a]{$3$}\taput{\small $$}} {
\Tdot~[tnpos=r]{$\textcolor[rgb]{0.00,0.00,1.00}{4}
[1,3]$}\taput{\small $$} \pstree {\Tdot~[tnpos=a]{$1$}\tbput{\small
$$}} { \Tdot~[tnpos=r]{$\textcolor[rgb]{0.00,0.00,1.00}{5} [1,3,1]$}\taput{\small $$} } } \pstree
{\Tdot~[tnpos=a]{$4$}\taput{\small $$}} {
\Tdot~[tnpos=r]{$\textcolor[rgb]{0.00,0.00,1.00}{5}
[1,4]$}\taput{\small $$} } } }

 \caption{Branch of the tree for the construction of
the gop with $\omega_1=1$ on $\mathcal{G}(\mathcal{F}_5)$}
\label{fig3}
\end{figure}

For example, for $N=5$, we construct one branch of a tree with
$\omega_1=1$ (see Fig. \ref{fig3}) : each vertex is an ordered
orbit,
the modulus of the gop is written on the last edge.\\
However, the sequence of ordered gop differs from the natural
downward lecture of the tree and has to be done following the
algorithm.

\section{Cardinal of classes}\label{section4}

In this section we emphasize some closed formulas giving the
cardinal of classes of gop. Recalling first the already known
formula for the class $[\overline{1_{\widetilde{k}}}]_N$ for which
we give a detailed proof, we consider the case were the class
possesses exactly one $k$-cycle, the case with only two cycles
belonging to the class and finally the main general formula of any
cycles with any length. We give rigorous proof of all. The general
formula is very interesting in the sense that even using computer
network it is impossible to check every function of  $\mathcal{F}_N$
when $N$ is larger than 100.\\

\subsection{Discrete maps with 1-cycle only}

The theorem 4 gives the number of discrete  maps of $\mathcal{F}_N$
which have only fixed points and no cycles of
length greater than one. This formula is explicit in \cite{Knuth} and \cite{Purd}. A complete proof is given here in detail.\\

 \label{th4}\thm{ Let $k$ be an integer between 1 and $N$. The number of
functions whose global orbit pattern is $ [1_{\widetilde{k}}]_N$
(i.e. belonging to the class $[\overline{1_{\widetilde{k}}}]_N$) is
$\left(\begin{array}{c}
  N-1 \\
  N-k \\
\end{array}\right) N^{N-k}$.}

That is \begin{equation}\label{eq11}
 \sharp [\overline{1_{\widetilde{k}}}]_N = \left(\begin{array}{c}
  N-1 \\
  N-k \\
\end{array}\right) N^{N-k}.
\end{equation}

\proof{ Let $k$ be a non-zero integer. Let $f$ be a function of
$\mathcal{F}_N$. There are $\left(\begin{array}{c}
  N \\
  k \\
\end{array}\right)$ possibilities to choose $k$ fixed points. There remain
$N-k$ points. Let $p$ be an integer between $1$ and $N-k$. We assume
that $p$ points are directly connected to the $k$ fixed points. For
each of them, there are $k$ manners to choose one fixed point. There
are $k^p$ ways to connect directly $p$ points to $k$ fixed points.
There remains $N-k-p$ points that we must connect to the $p$ points.
There are $\sharp [1_{\widetilde{p}}]_{N-k}$ functions. Finally, the
number of functions with $k$ fixed points is $\displaystyle {N
\choose k} \sum\limits_{p=1}^{N-k} k^p ~\sharp
[\overline{1_{\widetilde{p}}}]_{N-k}$. We now prove recursively on
$N$ for every $0 \leq k \leq N$ that $\sharp
[\overline{1_{\widetilde{k}}}]_N = \left(\begin{array}{c}
  N-1 \\
  N-k \\
\end{array}\right) N^{N-k}$.
We have $\sharp \overline{[1]}_1=1$. The formula is true.\\ We
suppose that $\forall k \leq N$ $\sharp
[\overline{1_{\widetilde{k}}}]_N = \left(\begin{array}{c}
  N-1 \\
  N-k \\
\end{array}\right) N^{N-k}$.\\
Let $X$ be a set with $N+1$ elements. We look for the functions of
$\mathcal{F}_{N+1}$ which have $k$ fixed points. Thanks to the
previous reasoning, we have \\
$\sharp [\overline{1_{\widetilde{k}}}]_{N+1} = \displaystyle
\left(\begin{array}{c}
  N+1 \\
  k \\
\end{array}\right)
\sum\limits_{p=1}^{N+1-k} k^p ~\sharp [\overline{1_{\widetilde{p}}}]_{N+1-k}$.\\
$\sharp [\overline{1_{\widetilde{k}}}]_{N+1} = \displaystyle
\left(\begin{array}{c}
  N+1 \\
  k \\
\end{array}\right)
\sum\limits_{p=1}^{N+1-k} k^p ~\sharp [\overline{1_{\widetilde{p}}}]_{N-(k-1)}$.\\
We use the recursion assumption.\\
$\sharp [\overline{1_{\widetilde{k}}}]_{N+1} = \displaystyle {N+1
\choose
k} \sum\limits_{p=1}^{N+1-k} k^p {N-k \choose p-1} (N-k+1)^{N-k+1-p}$.\\
$\sharp [\overline{1_{\widetilde{k}}}]_{N+1} =\displaystyle k {N+1
\choose
k} \sum\limits_{p=0}^{N-k}{N-k \choose p}~ k^{p} (N-k+1)^{N-k-p}$.\\
$\sharp [\overline{1_{\widetilde{k}}}]_{N+1} =\displaystyle k {N+1
\choose
k} (N+1)^{N-k}$.\\
$\sharp [\overline{1_{\widetilde{k}}}]_{N+1} =\displaystyle {N
\choose k-1}
(N+1)^{N-k+1}$.\\
$\sharp [\overline{1_{\widetilde{k}}}]_{N+1} =\displaystyle {N
\choose N-k+1} (N+1)^{N-k+1}$. q.e.d.}

\subsection{Discrete maps with $k$-cycle}
We look now for the number of functions with exactly one $k$-cycle.\\

\thm{Let $k$ be an integer between 1 and $N$. The number of
functions whose global orbit pattern is $[k]_N$ is $\sharp
[\overline{1_{\widetilde{k}}}]_{N} \times (k-1)!$.}

i.e. \begin{equation}\label{eq12} \sharp \overline{[k]}_N = \sharp
[\overline{1_{\widetilde{k}}}]_{N} \times (k-1)!.
\end{equation}

\proof{ There are ${N \choose k}$ ways of choosing $k$ elements
among $N$. Then, there are $(k-1)!$ choices for the image of those
$k$ elements in order to constitute a $k$-cycle by $f$. We must now
count the number of ways of connecting directly or not the remaining
$N-k$ elements to the $k$-cycle. We established already this number
which is equal to $\displaystyle \sum\limits_{p=1}^{N-k} k^p ~\sharp
[\overline{1_{\widetilde{p}}}]_{N-k}$. Finally, we have $\sharp
\overline{[k]}_N =(k-1)! {N \choose k} \sum\limits_{p=1}^{N-k} k^p
~\sharp [\overline{1_{\widetilde{p}}}]_{N-k}$. That is, $\sharp
[k]_N =\sharp [\overline{1_{\widetilde{k}}}]_{N} \times (k-1)!$.
q.e.d. }

\subsection{Discrete maps with only two cycles}

We give the number of functions with only two cycles.\\

\thm{Let $N \geq 2$. Let $p$ and $q$ be two non-zero integers such
that $p+q \leq N$. Then,
\begin{equation}\label{eq13}
 \sharp[\overline{ p,q}]_N = \sharp [\overline{1_{\widetilde{p+q}}}]_N \frac{(p+q-1)!}{q}= \frac{(N-1)!~
N^{N-(p+q)}}{(N-(p+q))!~ q}.
\end{equation}}\\

\proof{ We consider a function $f$ which belongs to the class
$[\overline{1_{\widetilde{p+q}}}]_N$. We search the number of
functions constructed from $f$ whose gop is $[p,q]_N$. From the $p$
fixed points of $f$, we construct a $p$-cycle. Thus, there are
${{p+q-1} \choose {p-1}}$ ways to choose $p-1$ integers among the
$p+q-1$ fixed points. Counting the first given fixed point of $f$,
we have $p$ points which allow to construct $(p-1)!$ functions with
a $p$-cycle. Then there remain $q$ points which give $(q-1)!$
different functions with a $q$-cycle. Finally, the number of
functions whose gop is $[p,q]_N$ is : ${{p+q-1} \choose {p-1}}
(p-1)! (q-1)!$ that is the formula
$\frac{(p+q-1)!}{q}$. }\\

\rem{We notice that for all $k$ non-zero integer such that $ k \leq
N-1$, ~~ $\sharp \overline{[k,1]}_N = \sharp \overline{[k+1]}_N$.}

\subsection{General case : discrete  maps with cycles of any length}

We introduce now the main theorem of the section which gives the
number of gop of discrete maps thanks to a closed formula.\\

Given a global orbit pattern $\alpha$, the next theorem gives a
formula which gives
the number of functions which belong to $\overline{\alpha}$.\\

\thm{Let $p \geq 2$ be an integer. Let $ [\omega_1, \ldots,
\omega_p]_N$ be a gop of $ \mathcal{G}(\mathcal{F}_N)$. Then, \\
\begin{equation}\label{eq14}
\sharp \overline{[\omega_1, \ldots, \omega_p]}_N = \sharp
[\overline{1_{\widetilde{\omega_1 + \ldots +\omega_p}} }]_N
\frac{(\omega_1 + \ldots +\omega_p -1)!} {\omega_p \times
(\omega_{p-1}+\omega_{p})\times \ldots \times (\omega_2 + \ldots
+\omega_p)}
\end{equation}
\begin{equation}\label{eq15}
 \sharp
\overline{[\omega_1, \ldots, \omega_p]}_N = \frac{(N-1)!
~N^{N-(\omega_1 + \ldots +\omega_p)}}{(N-(\omega_1 + \ldots
+\omega_p))! ~\prod\limits_{k=2}^{p} (\sum\limits_{j=k}^{p}
\omega_j)}
\end{equation}}

\proof{ We consider a function $f$ which belongs to $
[\overline{1_{\widetilde{\omega_1 + \ldots +\omega_p}}}]_N$. We
search the number of functions constructed from $f$ whose gop is
$[\omega_1, \ldots, \omega_p]_N $. From the $\omega_1$ fixed points
of $f$, we construct a $\omega_1$-cycle. Thus, there are ${{\omega_1
+ \ldots +\omega_p-1} \choose {\omega_1-1}}$ ways to choose
$\omega_1-1$ integers among the $\omega_1 + \ldots +\omega_p-1$
fixed points. Counting the first given fixed point of $f$, we have
$\omega_1$ points which allow to construct $(\omega_1-1)!$ functions
with a $\omega_1$-cycle. Then, the first fixed point of $f$ which
has not be chosen for the $\omega_1$-cycle, will belong to the
$\omega_2$-cycle. Thus, there are ${{\omega_2 + \ldots +\omega_p-1}
\choose {\omega_2-1}}$ ways to choose $\omega_2-1$ integers among
the $\omega_2 + \ldots +\omega_p-1$ fixed points. So we have
$\omega_2$ points which allow to construct $(\omega_2-1)!$ functions
with a $\omega_2$-cycle. We keep going on that way until there
remain $\omega_p$ fixed points which allow to construct
$(\omega_p-1)!$ functions with a $\omega_p$-cycle. Finally, we have
constructed : \\${{\omega_1 + \ldots +\omega_p-1} \choose
{\omega_1-1}} (\omega_1-1)! {{\omega_2 + \ldots +\omega_p-1} \choose
{\omega_2-1}} (\omega_2-1)! \times \ldots \times
{{\omega_{p-1}+\omega_p-1} \choose {\omega_{p-1}-1}}
(\omega_{p-1}-1)!(\omega_p-1)!$ functions. We simplify and obtain
the formula. }\\

\cor{Let $p$ be a non-zero integer. Let $ [\omega_1, \ldots,
\omega_p]_N$ be a gop of $\mathcal{G}(\mathcal{F}_N)$. We suppose
that there exists $j$ such that $\omega_j \geq 2$. Let $h$ be an
integer between 1 and $\omega_j-1$. Then
\begin{equation}\label{eq16}
\sharp \overline{[\omega_1, \ldots, \omega_j, \ldots, \omega_p]}_N =
\sharp \overline{[\omega_1, \ldots, \omega_j-h,h, \omega_{j+1},
\ldots, \omega_p]}_N \times ( h + \omega_{j+1} + \ldots +
\omega_p).
\end{equation}}\\

\proof{ $\sharp \overline{[\omega_1, \ldots, \omega_j-h,h,
\omega_{j+1}, \ldots, \omega_p]}_N \times ( h + \omega_{j+1} +
\ldots + \omega_p) = \sharp[ \overline{1_{\widetilde{\small{\omega_1
+ \ldots +\omega_p}}}}]_N
\\ \times \frac{(\omega_1 + \ldots +\omega_p -1)! ( h + \omega_{j+1} +
\ldots + \omega_p)}{\omega_p (\omega_{p-1}+\omega_{p}) \ldots
(\omega_{j+1}+\ldots + \omega_{p}) (h + \omega_{j+1}+\ldots +
\omega_{p}) (\omega_j + \omega_{j+1}+\ldots + \omega_{p}) \times
\ldots \times (\omega_2 +
\ldots +\omega_p)}$. \\

We simplify and we exactly obtain \\

$\sharp \overline{[\omega_1, \ldots, \omega_j-h,h, \omega_{j+1},
\ldots, \omega_p]}_N \times ( h + \omega_{j+1} + \ldots + \omega_p)
= \sharp \overline{[\omega_1,
\ldots, \omega_j, \ldots, \omega_p]}_N$. }\\

Examples :\\

$\sharp \overline{[2_{\widetilde{2}},1,3]}_{11} = 11,180,400$.\\

$\sharp \overline{[5,2,10,8,15,2,3]}_{50} = 29, 775, 702, 147, 667,
389, 218, 762, 343, 520, 975, 006, \\348,329, 578, 044, 480, 000,
000, 000,
000, 000$.\\

$\sharp \overline{[5,2,10,8,15,2,3]}_{50} \cong 2.98 \times 10^{63}$
among the $8.88 \times 10^{84}$ functions of $\mathcal{F}_{50}$.

\section{Functions with local properties }\label{section5}

\subsection{Locally rigid functions}

Obviously it is not possible to transpose to the functions on finite
sets the notions of continuity and derivability which play a
dramatic role in mathematical analysis since several centuries. In
fact the class $\mathcal{C}_0(I)$ of the continuous functions on the
real interval $I$ is a very small subset of the set $I^{\mathbb{R}}$
of all the functions on $I$. Hence by analogy to this fact and
trying to mimic some others properties of continuous functions, we
introduce some subsets of particular functions of $\mathcal{F}_{N}$,
which have local properties such as locally bounded range in a sense
we precise further. Limiting the range of the function in a
neighbourhood of any point of the interval induces a kind of
"rigidity" of the function, hence we call these functions locally
rigid functions. In these subsets, the gop are found to be fully
efficient in order to describe very precisely the dynamics of the
orbits. We first consider the very simple subset
$\mathcal{L}\scriptstyle{\mathcal{R}}_{1,N}$ of functions for which
the difference between $f(p)$ and $f(p+1)$ is drastically bounded.
In next subsection we consider more sophisticated
subsets.\\

We consider the set : \\

$\mathcal{L}\scriptstyle{\mathcal{R}}_{1,N} = \{ f \in
\mathcal{F}_N$ such that $ \forall p, 0
\leq p \leq N-2, |f(p)-f(p+1)| \leq 1 \}$.\\

\subsubsection{Orbits of $\mathcal{L}\scriptstyle{\mathcal{R}}_{1,N}$} \label{sec:52}

\thm{If $f \in \mathcal{L}\scriptstyle{\mathcal{R}}_{1,N}$ then $f$
has only periodic orbits of order $1$ or $2$.}

\proof{We suppose that $f \in
\mathcal{L}\scriptstyle{\mathcal{R}}_{1,N}$ has a 3-cycle. We denote
$(a; f(a);f^2(a))$ taking $a$ the smallest value of the 3-cycle. If
$a < f(a)<f^2(a)$ then there exist two non-zero integers $e$ and
$e'$ such that $f(a)=a+e$ and $f^2(a)=f(a)+e'$. Thus, $f^2(a)-e'
\leq f^3(a) \leq f^2(a)+e'$. That is $f(a)\leq a \leq f(a)+2e'$. And
finally we have the relation $a +e \leq a$ which is
impossible.\\
If $a < f^2(a)<f(a)$ then there exist two non-zero integers $e$ and
$e'$ such that $f^2(a)=a+e$ and $f(a)=f^2(a)+e'$. Thus, $f(a)-e \leq
f^3(a) \leq f(a)+e$. That is $f(a)-e\leq a \leq f(a)+e$. But
$f(a)-e=a+e'$. And finally we have the relation $a +e' \leq a$ which
is impossible.\\
We can prove in the same way that the function $f$ can't have either
3-cycle or greater order cycle than 3.}

\subsubsection{Numerical results and conjectures }
We have done numerical studies of the
$\mathcal{G}(\mathcal{L}\scriptstyle{\mathcal{R}}_{1,N})$ for $N =
1$ to $16$, using the brute force of a desktop computer (i.e.
checking every function
belonging to these sets).\\

The Tables \ref{tab:tab11}, \ref{tab:tab12}, \ref{tab:tab13},
\ref{tab:tab14}, \ref{tab:tab15} and \ref{tab:tab16} show the
sequences for
$\mathcal{L}\scriptstyle{\mathcal{R}}_{1,1}$ to $\mathcal{L}\scriptstyle{\mathcal{R}}_{1,16}$.\\

In theses Tables we display in the first column all the gop of
$\mathcal{G}(\mathcal{L}\scriptstyle{\mathcal{R}}_{1,N})$ for every
value of $N$. For a given $N$, there are two columns; the left one
displays the cardinal of every existing class of gop (- stands for
non existing gop). Instead the second shows more regularity,
displaying on the row of the gop $[2_{\widetilde{k}}]$ the sum of
the cardinals of all the classes of the gop of the form
$[\underbrace{2,2, \ldots
\ldots,\underbrace{1}_{i^{\small{\textrm{th}}}},\ldots, 2}_{k+1
~\textrm{orders}}]$ which exist.

Then we are able to formulate some statements which have not yet
been proved.\\

\begin{table}
\begin{tabular}{l*{8}{r}}
  \hline
g.o.p.&N=1&N=1&N=2&N=2&N=3&N=3&N=4&N=4 \\
Total number&&1&&4&&17&&68\\
$[1]$&1&+&2&+&7&+&26&+\\
$[1_{\widetilde{2}}]$&-&+&1&+&4&+&14&+\\
 $[1_{\widetilde{3}}]$&-&+&-&+&1&+&4&+\\
 $[1_{\widetilde{4}}]$&-&+&-&+&-&+&1&+\\
 $[ 2]$&-&+&1&1&4&4&18&18\\
$[2,1]$&-&+&-&+&1&1&3&4\\
$[1,2]$&-&+&-&+&-&+&1&+ \\
$[2_{\widetilde{2}}]$&-&+&-&+&-&+&1&1 \\
   \hline
\end{tabular}
\caption{Numbering the locally rigid functions for $f \in
\mathcal{L}\scriptstyle{\mathcal{R}}_{1,1}$, $f \in
\mathcal{L}\scriptstyle{\mathcal{R}}_{1,2}$, $f \in
\mathcal{L}\scriptstyle{\mathcal{R}}_{1,3}$, $f \in
\mathcal{L}\scriptstyle{\mathcal{R}}_{1,4}$. } \label{tab:tab11}
\end{table}

\begin{table}
\begin{tabular}{l*{6}{r}}
  \hline
g.o.p.&N=5&N=5&N=6&N=6&N=7&N=7 \\
Total number&&259&&950&&387\\
$[1]$&95&+&340&+&1,193&+ \\
$[1_{\widetilde{2}}]$&50&+&174&+&600&+\\
 $[1_{\widetilde{3}}]$&16&+&58&+&204&+ \\
 $[1_{\widetilde{4}}]$&4&+&16&+&60&+ \\
 $[1_{\widetilde{5}}]$&1&+&4&+&16&+\\
$[1_{\widetilde{6}}]$&-&+&1&+&4&+ \\
$[1_{\widetilde{7}}]$&-&+&-&+&1&+ \\
$[ 2]$&70&70&264&264&952&952\\
$[2,1]$&12&18&45&70&166&264 \\
$[1,2]$&6&+&25&+&98&+ \\
$[2_{\widetilde{2}}]$&4&4&18&18&70&70 \\
$[2_{\widetilde{2}},1]$&1&1&4&4&17&18\\
 $[1,2_{\widetilde{2}}]$&-&+&-&+&1& +\\
$[ 2,1,2]$&-&+&-&+&-&+\\
 $[ 2_{\widetilde{3}}]$&-&+&1&1&4&4 \\
  $[2_{\widetilde{3}},1]$&-&+&-&+&1&1\\

   \hline
\end{tabular}
\caption{Numbering the locally rigid functions for $f \in
\mathcal{L}\scriptstyle{\mathcal{R}}_{1,5}$, $f \in
\mathcal{L}\scriptstyle{\mathcal{R}}_{1,6}$, $f \in
\mathcal{L}\scriptstyle{\mathcal{R}}_{1,7}$. } \label{tab:tab12}
\end{table}

\begin{table}
\begin{tabular}{l*{6}{r}}
  \hline

g.o.p.&N=8&N=8&N=9&N=9&N=10&N=10 \\
Total number&&11,814&&40,503&&13,6946 \\
$[1]$&4,116&+&14,001&+&47,064&+\\
$[ 1_{\widetilde{2}}]$&2,038&+&6,852&+&22,806& +\\
$[ 1_{\widetilde{3}}]$&700&+&2,366&+&7,896&+\\
$[ 1_{\widetilde{4}}]$&214&+&742&+&2,520& +\\
$[1_{\widetilde{5}}]$&60&+&216&+&754& +\\
$[ 1_{\widetilde{6}}]$&16&+&60&+&216&+\\
$[1_{\widetilde{7}}]$&4&+&16&+&60& +\\
$[1_{\widetilde{8}}]$&1&+&4&+&16&+ \\
$[1_{\widetilde{9}}]$&-&+&1&+&4&+ \\
$[1_{\widetilde{10}}]$&-&+&-&+&1&+\\
$[ 2]$&3,356&3,356&11,580&11,580&39,364&39,364\\
$[2,1]$&590&952&2,062&3,356&7,072&11,580 \\
$[1,2]$&362&+&1,294&+&4,508&+\\
$[ 2_{\widetilde{2}}]$&264&264&952&952&3,356&3,356\\
$[2_{\widetilde{2}},1]$&62&70&222&264&770&952 \\
$[1,2_{\widetilde{2}}]$&6&+&28&+&113&+ \\
$[2,1,2]$&2&+&14&+&69&+\\
$[ 2_{\widetilde{3}}]$&18&18&70&70&264&264 \\
 $[2_{\widetilde{3}},1]$&4&4&18&18&69&70\\
$[1,2_{\widetilde{3}}]$&-&+&-&+&1& +\\
$[2,1,2_{\widetilde{2}}]$&-&+&-&+&-&+ \\
$[2_{\widetilde{2}},1,2]$&-&+&-&+&-&+ \\
$[2_{\widetilde{4}}]$&1&1&4&4&18&18 \\
$[2_{\widetilde{4}},1]$&-&+&1&1&4&4 \\
$[1,2_{\widetilde{4}}]$&-&+&-&+&-& +\\
$[2,1,2_{\widetilde{3}}]$&-&+&-&+&-& +\\
$[2_{\widetilde{2}},1,2_{\widetilde{2}}]$&-&+&-&+&-& +\\
$[2_{\widetilde{3}},1,2]$&-&+&-&+&-& +\\
$[2_{\widetilde{5}}]$&-&+&-&+&1&1 \\
 \hline
\end{tabular}
\caption{Numbering the locally rigid functions for $f \in
\mathcal{L}\scriptstyle{\mathcal{R}}_{1,8}$, $f \in
\mathcal{L}\scriptstyle{\mathcal{R}}_{1,9}$, $f \in
\mathcal{L}\scriptstyle{\mathcal{R}}_{1,10}$. } \label{tab:tab13}
\end{table}

\begin{table}
\begin{tabular}{l*{6}{r}}
  \hline
g.o.p.&N=11&N=11&N=12&N=12&N=13&N=13 \\
Total number&&457,795&&1,515,926&&4,979,777 \\
$[1]$&156,629&+&516,844&+&1,693,073&+\\
$[1_{\widetilde{2}}]$&75,292&+&246,762&+&803,706&+\\
$[1_{\widetilde{3}}]$&26,098&+&85,556&+&278,580& +\\
$[1_{\widetilde{4}}]$&8,434&+&27,904&+&91,488&+\\
$[1_{\widetilde{5}}]$&2,756&+&8,658&+&28,738&+ \\
$[1_{\widetilde{6}}]$&756&+&2,590&+&8,730&+\\
$[1_{\widetilde{7}}]$&216&+&756&+&2,592&+ \\
$[1_{\widetilde{8}}]$&60&+&216&+&756&+\\
$[1_{\widetilde{9}}]$&16&+&60&+&216&+ \\
$[1_{\widetilde{10}}]$&4&+&16&+&60&+ \\
$[1_{\widetilde{11}}]$&1&+&4&+&16&+ \\
$[1_{\widetilde{12}}]$&-&+&1&+&4& +\\
$[1_{\widetilde{13}}]$&-&+&-&+&1&+\\
$[2]$&132,104&132,104&438,846&438,846&1,445,258&1,445,258\\
$[2,1]$&23,941&39,364&80,108&132,104&265,548&438,846 \\
$[1,2]$&15,423&+&51,996&+&173,298&+\\
$[2_{\widetilde{2}}]$&11,580&11,580&39,364&39,364&132,104&132,104\\
$[2_{\widetilde{2}},1]$&2,634&3,356&8,883&11,580&29,659&39,364 \\
$[1,2_{\widetilde{2}}]$&429&+&1,555&+&5,478&+\\
$[2,1,2]$&293&+&1,142&+&4,227& +\\
$[2_{\widetilde{3}}]$&952&952&3,356&3,356&11,580&11,580\\
$[2_{\widetilde{3}},1]$&255&264&899&952&3,098&3,356\\
$[ 1,2_{\widetilde{3}}]$&7&+&35&+&152&+ \\
$[ 2,1,2_{\widetilde{2}}]$&2&+&16&+&86&+\\
$[2_{\widetilde{2}},1,2]$&-&+&2&+&20& +\\
$[2_{\widetilde{4}}]$&70&70&264&264&952&952\\
$[ 2_{\widetilde{4}},1]$&18&18&70&70&263&264 \\
$[ 1,2_{\widetilde{4}}]$&-&+&-&+&1&+\\
$[2,1,2_{\widetilde{3}}]$&-&+&-&+&-&+ \\
$[2_{\widetilde{2}},1,2_{\widetilde{2}}]$&-&+&-&+&-&+ \\
$[2_{\widetilde{3}},1,2]$&-&+&-&+&-&+ \\
$[2_{\widetilde{5}}]$&4&4&18&18&70&70 \\
$[2_{\widetilde{5}},1]$&1&1&4&4&18&18 \\
$[1,2_{\widetilde{5}}]$&-&+&-&+&-&+ \\
$[2_{\widetilde{6}}]$&-&+&1&1&4&4 \\
$[2_{\widetilde{6}},1]$&-&+&-&+&1&1\\
 \hline
\end{tabular}
\caption{Numbering the locally rigid functions for $f \in
\mathcal{L}\scriptstyle{\mathcal{R}}_{1,11}$, $f \in
\mathcal{L}\scriptstyle{\mathcal{R}}_{1,12}$, $f \in
\mathcal{L}\scriptstyle{\mathcal{R}}_{1,13}$.} \label{tab:tab14}
\end{table}

\begin{table}
\begin{tabular}{l*{4}{r}}
  \hline
g.o.p.&N=14&N=14&N=15&N=15 \\
Total number&&16,246,924&&52,694,573
\\
$[1]$&5,511,218&+&17,841,247& +\\
$[1_{\widetilde{2}}]$&2,603,258&+&8,391,360&+\\
$[ 1_{\widetilde{3}}]$&901,802&+&2,904,592&+\\
$[1_{\widetilde{4}}]$&297,728&+&962,888& +\\
$[1_{\widetilde{5}}]$&94,440&+&307,848&+\\
$[ 1_{\widetilde{6}}]$&29,050&+&95,676&+\\
$[1_{\widetilde{7}}]$&8,746&+&29,140&+ \\
$[1_{\widetilde{8}}]$&2,592&+&8,748&+\\
$[ 1_{\widetilde{9}}]$&756&+&2,592&+\\
$[1_{\widetilde{10}}]$&216&+&756&+ \\
$[1_{\widetilde{11}}]$&60&+&216&+\\
$[ 1_{\widetilde{12}}]$&16&+&60&+\\
$[ 1_{\widetilde{13}}]$&4&+&16&+ \\
$[ 1_{\widetilde{14}}]$&1&+&4&+\\
$[1_{\widetilde{15}}]$&-&+&1& +\\
$[2]$&4,725,220&4,725,220&15,352,392&15,352,392\\
$[2,1]$&873,149&1,445,258&2,851,350&+\\
$[1,2]$&572,109&+&1,873,870&+\\
$[2_{\widetilde{2}}]$&438,846&438,846&1,445,258&1,445,258\\
$[2_{\widetilde{2}},1]$&98,135&132,104&322,310&438,846\\
$[1,2_{\widetilde{2}}]$&18,873&+&63,967&+ \\
$[2,1,2]$&15,096&+&52,569&+\\
$[2_{\widetilde{3}}]$&39,364&39,364&132,104&132,104\\
$[2_{\widetilde{3}},1]$&10,460&11,580&34,845&39,364 \\
$[1,2_{\widetilde{3}}]$&605&+&2,282&+\\
$[2,1,2_{\widetilde{2}}]$&389&+&1,596&+\\
$[ 2_{\widetilde{2}},1,2]$&126&+&641&+\\
$[2_{\widetilde{4}}]$&3,356&3,356&11,580&11,580\\
$[2_{\widetilde{4}},1]$&942&952&3,292&3,356 \\
$[1,2_{\widetilde{4}}]$&8&+&44&+\\
$[2,1,2_{\widetilde{3}}]$&2&+&18&+\\
$[2_{\widetilde{2}},1,2_{\widetilde{2}}]$&-&+&2&+\\
$[2_{\widetilde{3}},1,2]$&-&+&-&+ \\
$[2_{\widetilde{5}}]$&264&264&952&952\\
$[2_{\widetilde{5}},1]$&70&70&264&264 \\
$[1,2_{\widetilde{5}}]$&-&+&-&+\\
$[2_{\widetilde{6}}]$&18&18&70&70\\
$[2_{\widetilde{6}},1]$&4&4&18&18\\
$[2_{\widetilde{7}}]$&1&1&4&4 \\
$[2_{\widetilde{7}},1]$&-&+&1&1 \\
      \hline
\end{tabular}
\caption{Numbering the locally rigid functions for  $f \in
\mathcal{L}\scriptstyle{\mathcal{R}}_{1,14}$, $f \in
\mathcal{L}\scriptstyle{\mathcal{R}}_{1,15}$. } \label{tab:tab15}
\end{table}

\label{stat1} \stat{
\begin{equation}\label{eq17}
\sharp
[\overline{1_{\widetilde{k}}}]_{\mathcal{L}\scriptstyle{\mathcal{R}}_{1,N}}
=\sharp
[\overline{1_{\widetilde{k+1}}}]_{\mathcal{L}\scriptstyle{\mathcal{R}}_{1,N+1}}
~ \textrm{for} ~k \leq \frac{N+1}{2}.
\end{equation}}\\

\label{stat2} \stat{
\begin{equation}\label{eq18}
\sharp
[\overline{2_{\widetilde{k}}}]_{\mathcal{L}\scriptstyle{\mathcal{R}}_{1,N}}
= \sharp
[\overline{2_{\widetilde{k+1}}}]_{\mathcal{L}\scriptstyle{\mathcal{R}}_{1,N+2}}
~ \textrm{for} ~k \leq \frac{N}{2}.
\end{equation}}\\

 \label{stat3} \stat{ \begin{equation}\label{eq19}
 \sharp
[\overline{2_{\widetilde{k}}}]_{\mathcal{L}\scriptstyle{\mathcal{R}}_{1,N}}
=\sharp
[\overline{2_{\widetilde{k}},1}]_{\mathcal{L}\scriptstyle{\mathcal{R}}_{1,N+1}}~
\textrm{for} ~2k \leq N \leq 3k-1.
\end{equation}}\\

\label{stat4} \stat{
\begin{eqnarray}
  \sharp
[\overline{2_{\widetilde{k}}}]_{\mathcal{L}\scriptstyle{\mathcal{R}}_{1,N}}
 &=& \sum\limits_{i=1}^{k+1} \sharp [\overline{\underbrace{2,2, \ldots
\ldots,\underbrace{1}_{i^{\small{\textrm{th}}}},\ldots, 2}_{k+1
~\textrm{orders}}}]_{\mathcal{L}\scriptstyle{\mathcal{R}}_{1,N}} ~\textrm{for}~ 2k+1 \leq N \nonumber \\
  &=& \sum\limits_{i=1}^{k+1} \sharp
[\overline{2_{\widetilde{i-1}},1,2_{\widetilde{k-i+1}}}]_{\mathcal{L}\scriptstyle{\mathcal{R}}_{1,N}}
~\textrm{for}~ 2k+1 \leq N
\end{eqnarray}
}\\

\label{stat5} \stat{\begin{equation}\label{eq20}\sharp
[\overline{1_{\widetilde{N-k+1}}}]_{\mathcal{L}\scriptstyle{\mathcal{R}}_{1,N}}
=\left\{
                                          \begin{array}{ll}
                                            1 & {\textrm{if}}~~ k=1 \\
                                            2 & {\textrm{if}}~~ k=2 \\
                                            \left( \frac{4}{27}\right)(k+1) \times 3^k
& {\textrm{for}}~~ 3\leq k \leq \frac{N+1}{2}
\end{array}
                                        \right.
                                        \end{equation}}\\

\rem{We call $u_k= \sharp
[\overline{1_{\widetilde{N-k+1}}}]_{\mathcal{L}\scriptstyle{\mathcal{R}}_{1,N}}$.
For $k>2$, then $u_k$ is the sequence A120926 On-line Encyclopedia
of integer Sequences : it is the number of sequences where 0 is
isolated in
ternary words of length $N$ written with $\{0,1,2\}$. }\\

These statements show that first the set
$\mathcal{L}\scriptstyle{\mathcal{R}}_{1,N}$ is an interesting set
to be considered for dynamical systems and secondly
the gop are fruitful in this study. However the set\\

$\mathcal{L}\scriptstyle{\mathcal{R}}_{2,N} = \{ f \in
\mathcal{F}_N$ such that $ \forall p,
0 \leq p \leq N-2, |f(p)-f(p+1)| \leq 2 \}$\\

\noindent is too much large to give comparable results. Then we
introduce more sophisticated sets we call sets with locally bounded
range which more or less correspond to an analogue of the discrete
convolution product of the local variation of $f$ with a compact
support function $\overrightarrow{\alpha_t}$.

\begin{table}
\begin{tabular}{l*{2}{r}}
  \hline
g.o.p.&N=16&N=16 \\
Total number& & 170,028,792\\
$[1]$& 57,477,542&+\\
$[1_{\widetilde{2}}]$& 26,932,398&+\\
$[ 1_{\widetilde{3}}]$& 9,314,088&+\\
$[1_{\widetilde{4}}]$&3,097,650&+ \\
$[1_{\widetilde{5}}]$&996,764 &+\\
$[ 1_{\widetilde{6}}]$& 312,456&+\\
$[1_{\widetilde{7}}]$& 96,096&+\\
$[1_{\widetilde{8}}]$&29,158&+ \\
$[ 1_{\widetilde{9}}]$&8,748&+\\
$[1_{\widetilde{10}}]$& 2,592&+\\
$[1_{\widetilde{11}}]$& 756&+\\
$[ 1_{\widetilde{12}}]$& 216&+\\
$[ 1_{\widetilde{13}}]$&  60&+ \\
$[ 1_{\widetilde{14}}]$&16&+\\
$[1_{\widetilde{15}}]$& 4&+\\
$[1_{\widetilde{16}}]$& 1&+\\
$[2]$& 49,610,818 &49,610,818\\
$[2,1]$&9,255,822& 15,352,392\\
$[1,2]$& 6,096,570&+\\
$[2_{\widetilde{2}}]$&  4,725,220& 4,725,220\\
$[2_{\widetilde{2}},1]$&1,051,686& 1,445,258\\
$[1,2_{\widetilde{2}}]$&  213,975&+\\
$[2,1,2]$& 179,597&+\\
$[2_{\widetilde{3}}]$&438,846& 438,846 \\
$[2_{\widetilde{3}},1]$& 114,798& 132,104\\
$[1,2_{\widetilde{3}}]$& 8,284&+\\
$[2,1,2_{\widetilde{2}}]$& 6,146&+\\
$[ 2_{\widetilde{2}},1,2]$& 2,876 &+\\
$[2_{\widetilde{4}}]$&39,364& 39,364\\
$[2_{\widetilde{4}},1]$& 11,246& 11,580 \\
$[1,2_{\widetilde{4}}]$&204&+ \\
$[2,1,2_{\widetilde{3}}]$&106&+\\
$[2_{\widetilde{2}},1,2_{\widetilde{2}}]$& 22&+\\
$[2_{\widetilde{3}},1,2]$& 2&+\\
$[2_{\widetilde{5}}]$&3,356&3,356 \\
$[2_{\widetilde{5}},1]$&951&952\\
$[1,2_{\widetilde{5}}]$&1&+\\
$[2_{\widetilde{6}}]$&   264&264 \\
$[2_{\widetilde{6}},1]$&  70&70\\
$[2_{\widetilde{7}}]$&   18 & 18\\
$[2_{\widetilde{7}},1]$& 4&4 \\
$[2_{\widetilde{8}}]$& 1&1\\
      \hline
\end{tabular}
\caption{Numbering the locally rigid functions for $f \in
\mathcal{L}\scriptstyle{\mathcal{R}}_{1,16}$. } \label{tab:tab16}
\end{table}

\subsection{Orbits and patterns of locally rigid function sets}

Consider now the set :

 $\mathcal{L}\scriptstyle{\mathcal{R}}_{\overrightarrow{\alpha_t},q,N} = \{ f \in
\mathcal{F}_N$ such that $ \forall p, 0 \leq p \leq N-r-1,
\sum\limits_{r=1}^{r=t} \alpha_r |f(p)-f(p+r)| \leq q \} \bigcap \{
f \in \mathcal{F}_N$ such that $ \forall p, t \leq p \leq N-1,
\sum\limits_{r=1}^{r=t} \alpha_r |f(p)-f(p-r)| \leq q \}$ for the
vector $\overrightarrow{\alpha_t}=(\alpha_1, \alpha_2, \ldots,
\alpha_t) \in \mathbb{N}^t$, for $q \in \mathbb{N}$.

\begin{longtable}{*{5}{r}}

\caption{Numerical study of the set
$\mathcal{L}\scriptstyle{\mathcal{R}}_{\vec{\alpha_t},q,N}$ for
$N=10$, $t=5$, $\alpha_1=20$, $\alpha_2=9$, $\alpha_3=5$,
$\alpha_4=2$ and $\alpha_5=1$, for
$q=20, ..., 142$} \label{tab:tab17}\\
  \hline
q& maximal period &modulus& gop number& functions number \\
\hline
\endfirsthead
\caption[]{(Next)}\\
\hline
q& maximal period &modulus& gop number& functions number \\
\hline
\endhead
 \multicolumn{5}{r}{Following next page}\\ \hline
\endfoot
\multicolumn{5}{r}{
}
\endlastfoot
20 & 1  & 1  & 1 &  10\\
26  &2  & 2&   3 &  82\\
44  &2 &  3 &  6  & 21,764 \\
49  &3& 3& 7 &  48,112\\
50 & 3& 3& 7 &  53,210 \\
56 & 3 &4& 9  & 208,692 \\
59  &4 &  4& 15 &330,800 \\
63 & 4 &  5 &19 &626,890 \\
66 & 4 &10 & 37& 952,228 \\
67&  4 &10& 46& 1,064,316 \\
72 & 5 &10& 50& 1,630,018 \\
74 & 6 &10& 60& 1,816,826 \\
76 & 6 &10& 61 &2,152,450 \\
77 & 6 &10& 88& 2,416,368 \\
78  &6 &10& 91& 2,762,434 \\
79 &6 &10 &97 &3,188,080 \\
80 & 6& 10& 99 &3,735,666 \\
84  &6& 10& 100 &5,876,324 \\
85 & 6& 10& 103& 6,473,288 \\
87& 6 &10& 105& 7,851,728 \\
88 &7 &10& 121& 8,644,178 \\
89 &8& 10& 129 &9,521,920 \\
91 & 8 &10& 136& 11,414,556 \\
92 &8& 10& 165& 12,454,440 \\
94 &8 &10& 175 &14,756,058 \\
95& 8& 10& 177& 16,077,780 \\
96& 8& 10& 184 &17,208,654 \\
97& 8 &10& 185& 18,369,854 \\
98 &8 &10& 188 &19,585,746 \\
100 &8 &10& 192 &22,083,852 \\
101 &8 &10& 199 &23,584,452 \\
102 &8& 10& 204& 25,513,892 \\
103 &8 &10& 244& 27,912,772 \\
104 &8& 10& 304 &30,560,238 \\
105 &9& 10& 333 &33,516,466 \\
106 &9& 10& 380& 36,682,960 \\
107 &9 &10& 424 &40,004,280 \\
108& 10& 10& 491 &43,685,352 \\
109& 10& 10& 517& 47,655,856 \\
110& 10& 10& 529 &51,785,410 \\
111 & 10& 10& 562 &55,907,120 \\
112& 10& 10& 583& 60,341,276 \\
113& 10 &10& 612 &64,930,790 \\
114& 10 &10& 647& 69,766,178 \\
115& 10 &10& 706& 74,989,752 \\
116& 10 &10& 747& 80,087,120 \\
117& 10 &10& 791 &85,570,272 \\
118& 10 &10& 820& 91,206,218 \\
119& 10 &10& 836& 97,040,288 \\
120& 10 &10& 852& 103,121,916 \\
121& 10 &10& 872& 109,650,464 \\
122& 10 &10& 896& 116,345,296 \\
123& 10 &10& 919& 123,241,156 \\
124& 10 &10& 924& 130,360,938 \\
125& 10 &10& 928& 137,636,628 \\
126& 10 &10& 930& 145,536,068 \\
127& 10 &10& 932& 154,370,862 \\
128& 10 &10& 938& 164,145,928 \\
129& 10 &10& 960& 174,942,026 \\
130& 10 &10& 986 &186,438,038 \\
131& 10 &10& 1,006& 198,594,118 \\
132& 10 &10& 1,013& 211,550,402 \\
133& 10 &10& 1,015 &225,324,700 \\
134& 10 &10& 1,021 &239,976,118 \\
135& 10 &10&  1,022& 255,106,866 \\
137& 10 &10&  1,023& 286,726,234 \\
142& 10 &10&  1,023 &374,355,356\\
      \hline
\end{longtable}

The functions belonging to these sets show a kind of "rigidity": the
less is $q$, the more "rigid" is the function, this "rigidity" being
modulated by the vector $\overrightarrow{\alpha_t}$. Furthermore,
the maximal length of a periodic orbit increases with $q$, and so
the number of gop $\sharp
\mathcal{G}(\mathcal{L}\scriptstyle{\mathcal{R}}_{\overrightarrow{\alpha_t},q,N})$
and the maximal modulus of the gop.\\

\rem{Using this generalized notation, one has :
$\mathcal{L}{\scriptstyle{\mathcal{R}}}_{1,n}=\mathcal{L}{\scriptstyle{\mathcal{R}}}_{1,1,n}$
and
$\mathcal{L}{\scriptstyle{\mathcal{R}}}_{2,n}=\mathcal{L}{\scriptstyle{\mathcal{R}}}_{1,2,n}$.}

As an example, we explore numerically the case : $N=10$, $t=5$,
$\alpha_1=20$, $\alpha_2=9$, $\alpha_3=5$, $\alpha_4=2$ and
$\alpha_5=1$, for $q=20, \ldots, 142$. The results are displayed in
Table \ref{tab:tab17}. In this Table "modulus" means the maximal
modulus of the gop belonging to this set for the corresponding value
of $q$ in the row, "gop number" stands for $\sharp
\mathcal{G}(\mathcal{L}\scriptstyle{\mathcal{R}}_{\overrightarrow{\alpha_t},q,N})$
and "functions number" for $\sharp
\mathcal{L}\scriptstyle{\mathcal{R}}_{\overrightarrow{\alpha_t},q,N}$.
One can point out that for the particular function
$\overrightarrow{\alpha_t}$ of the example; it is possible to find
$10$ intervals $I_1, I_2, \ldots, I_{10} \subset \mathbb{N}$ such
that if $q \in I_r$ then there is no periodic orbit whose period is
strictly greater than $r$, (e.g., $I_6= \dbl 74,87 \dbr$).
Furthermore it is possible to split these intervals into
subintervals $I_{r,s}$ in which $\sharp
\mathcal{G}\left(\mathcal{L}\scriptstyle{\mathcal{R}}_{\overrightarrow{\alpha_t},q,N}\right)$
is constant when $q$ thumbs $I_{r,s}$. This is not the case for
$\sharp
\mathcal{L}\scriptstyle{\mathcal{R}}_{\overrightarrow{\alpha_t},q,N}$.

\section{Conclusion}\label{section6}
A discrete dynamical system associated to a function on finite
ordered set X can only exhibit periodic orbits. However the number
of the periods and the length of each are not easily predictable. We
formalise such a gop as the ordered set of periods when the initial
value thumbs $X$ in the increasing order. We can predict by means of
closed formulas, the number of gop of the set of all the function
from $X$ to itself. We also explore, using the brute force of
computers, some subsets of locally rigid functions on $X$, for which
interesting patterns of periodic orbits are found. Further study is
needed to understand the behaviour of dynamical systems associated
to functions belonging to these sets.



\newpage


\begin{thebibliography}{99}

\bibitem{Dia}
Diamond, P. and Pokrovskii, A., Statistical laws for computational
collapse of discretized chaotic mappings. \emph{International
Journal of Bifurcation and Chaos}, 1996, Vol.6, 12A, pp. 2389-2399.

\bibitem{Knuth}
Donald E. Knuth, \emph{The Art Of Computer Programming (Third
Edition)}, Addison-Wesley, 2005, Vol. 2, pp.8-9.



\bibitem{Krus}
Martin D. Kruskal, \emph{The expected number of components under a
random mapping function}, INST CNRS, June 1954, p.392.



\bibitem{Shark1}
A.N. Sharkovskii [1964], \emph{Coexistence of cycles of continuous
mapping of the line into itself. Ukrainian Math.}, 1995,
International Journal of Bifurcation and Chaos, No.5, pp.1263-1273.


\bibitem{Lozi1}
R. Lozi, \emph{The importance of strange attractors for industrial
mathematics}, Trends in Industrial and Applied Mathematics, 2002,
Proceedings of the 1st International Conference on Industrial and
Applied Mathematics of the Indian Subcontinent, A.H. Siddiqi and M.
Kocvara (Eds) Kluwer Academic Publisher, pp.275-303.



\bibitem{Lozi2}
R. Lozi, \emph{Giga-periodic Orbits for Weakly Coupled tent and
Logistic Discretized Maps}, Modern Mathematical Models, Methods and
Algorithms for Real World Systems, 2007, A.H. Siddiqi, I.S. Duff and
O. Christensen (Editors), Anamaya Publishers, New Delhi, India,
pp.81-124.


\bibitem{Lozi3}
R. Lozi, \emph{New Enhanced Chaotic Number Generators}, Indian
Journal of Industrial and Applied Mathematics, vol.1 No.1, pp.1-23.


\bibitem{Mutaf}
Ljuben R. Mutafchiev, \emph{Limit theorem concerning random mapping
patterns}, Combinatorica, 1988, Vol. 8, pp.345-356.



\bibitem{Purd}
Purdom, P.W, and Williams, J.H., \emph{Cycle length in a random
function}, Trans. of Amer. Math. Soc., Sept. 1968, Vol. 133,
 No.2, pp.547-551.


\bibitem{Patent2007}
(2007) US Patent 7,170,997 - Method of generating pseudo-random
numbers in an electronic device, and a method of encrypting and
decrypting electronic data.

\bibitem{Patent2006}
(2006) US Patent 6,999,445 - Multiple access communication system
using chaotic signals and method for generating and extracting
chaotic signal.

\bibitem{Patent2001}
(2001) U.S. Pat. No. 5,048,086 assigned to Hughes Aircraft Company
is related to an encryption system based on chaos theory. The system
uses the logistic equation $x_{n+1}=\mu x_n(1-x_n)$, which is a
mapping exhibiting chaos for certain values of $\mu$. In the
computations, floating-point operations are used.




\end{thebibliography}
\end{document}